\documentclass{article}

\usepackage{arxiv}

\usepackage[utf8]{inputenc} % allow utf-8 input
\usepackage[T1]{fontenc}    % use 8-bit T1 fonts
\usepackage{hyperref}       % hyperlinks
\usepackage{url}            % simple URL typesetting
\usepackage{booktabs}       % professional-quality tables
\usepackage{amsfonts}       % blackboard math symbols
\usepackage{nicefrac}       % compact symbols for 1/2, etc.
\usepackage{microtype}      % microtypography
\usepackage{doi}
\usepackage{enumitem}
\usepackage{style}
\usepackage{graphicx}
\usepackage{tikz}

\title{POLYDIM: A C++ library for POLYtopal DIscretization Methods}

\author{ \href{https://orcid.org/0000-0001-8642-4258}{\includegraphics[scale=0.06]{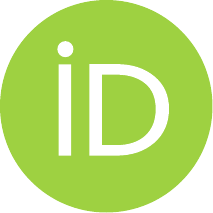}\hspace{1mm}Stefano~Berrone} \\
	Dipartimento di Scienze Matematiche\\
	``G. L. Lagrange''\\
	Politecnico di Torino, TO, 10129 \\
	\texttt{stefano.berrone@polito.it} \\
	\And
	\href{https://orcid.org/0000-0003-2016-5403}{\includegraphics[scale=0.06]{orcid.pdf}\hspace{1mm}Andrea~Borio} \\
	Dipartimento di Scienze Matematiche\\
	``G. L. Lagrange''\\
	Politecnico di Torino, TO, 10129 \\
	\texttt{andrea.borio@polito.it} \\
	\And
	\href{https://orcid.org/0000-0002-8540-3639}{\includegraphics[scale=0.06]{orcid.pdf}\hspace{1mm}Gioana~Teora} \\
	Dipartimento di Scienze Matematiche\\
	``G. L. Lagrange''\\
	Politecnico di Torino, TO, 10129 \\
	\texttt{gioana.teora@polito.it} \\
	\And
	\href{https://orcid.org/0000-0001-7123-9199}{\includegraphics[scale=0.06]{orcid.pdf}\hspace{1mm}Fabio~Vicini} \\
	Dipartimento di Scienze Matematiche\\
	``G. L. Lagrange''\\
	Politecnico di Torino, TO, 10129 \\
	\texttt{fabio.vicini@polito.it} \\
}

\renewcommand{\shorttitle}{POLYDIM}

\hypersetup{
pdftitle={POLYDIM},
pdfsubject={},
pdfauthor={Stefano~Berrone, Andrea~Borio, Gioana~Teora, Fabio~Vicini},
pdfkeywords={C++,implementation,Polytopal Method,Virtual Element Method},
}

\begin{document}
\maketitle

\begin{abstract}
This paper introduces PolyDiM, an open-source C++ library tailored for the development and implementation of polytopal discretization methods for partial differential equations. The library provides robust and modular tools to support advanced numerical techniques, with a focus on the Virtual Element Method in both 2D and 3D settings. PolyDiM is designed to address a wide range of challenging problems, including those involving non-convex geometries, Discrete Fracture Networks, and mixed-dimensional coupling. It is integrated with the geometry library GeDiM, and offers interfaces for MATLAB and Python to enhance accessibility. Distinguishing features include support for multiple polynomial bases, advanced stabilization strategies, and efficient local-to-global assembly procedures. PolyDiM aims to serve both as a research tool and a foundation for scalable scientific computing in complex geometrical settings.
\end{abstract}

% keywords can be removed
\keywords{C++ \and Implementation \and Polytopal Method \and Virtual Element Method}

\section{Introduction} \label{sec:introduction}

The numerical solution of partial differential equations (PDEs) on polytopal (polygonal or polyhedral) meshes has gained increasing attention in recent years due to its versatility and robustness in handling complex geometries. These meshes, which may include highly irregular or non-convex shapes, offer substantial advantages in engineering applications, since they can more naturally represent domains with intricate structures and reduce the computational complexities of some advanced meshing strategies such as refinement and agglomeration processes. Many novel numerical methods have been developed in recent years to deal with polytopal elements, such as the Virtual Element Method (VEM) \cite{DaVeigaBrezzi2023}, the Hybrid High-Order method (HHO) \cite{basicHHO}, and Polytopal Discontinuous Galerkin methods (PDG) \cite{PDG}.

In particular, the Virtual Element Method can be considered as a generalization of the Finite Element Method (FEM) which introduces in the local space suitable non-polynomial functions that are not known in a closed form, as well as the standard polynomials, by introducing projection and stability operators that make the discrete bilinear operators computable and stable. Among the advantages of VEM, in addition to the employment of polytopal meshes, we recall the possibility of handling discrete spaces of arbitrary $C^k$ global regularity \cite{Ck_regularity}, as well as spaces that satisfy exactly the divergence-free constraint \cite{Vacca2018}.

In this work, we present PolyDiM (POLYtopal DIscretization Method), an open-source C++ library tailored for the development and implementation of polytopal methods for PDEs. Built around the \textit{global-to-local} philosophy typical of the Finite Element Method (FEM) \cite{brennerscott}, PolyDiM is designed to be modular, extensible, and researcher-friendly. The library currently supports classical FEM in 1D, 2D, and 3D, alongside a robust implementation of VEM for a variety of boundary value problems in 2D and 3D.

PolyDiM is originally developed to support the research work of the authors and hope to be in support of the whole scientific community. Current applications include the simulation of coupled problems on intricate domains, such as Discrete Fracture Networks (DFNs) \cite{Vicini2024}, and 3D-1D coupling \cite{GrappeinTeora2025}. The library also supports research on mesh refinement and agglomeration strategies \cite{Sorgente2023}, alongside analysis of the robustness of VEM on badly-shaped elements \cite{Teora2024, Teora2024_mixed, Sorgente2022_2D, Sorgente2022_3D}.

To maximize interoperability and geometric flexibility, PolyDiM is fully integrated with GeDiM (GEometry for DIscretization MEthods) \cite{Gedim}, a C++ library developed by the same authors to provide robust geometry handling for polytopal domains. Furthermore, PolyDiM includes interfaces for MATLAB and Python, making it accessible to a broader range of users.

We recognize that several other software packages supporting polytopal discretization methods are currently available.
Table~\ref{tab:software} provides a summary of the most notable ones, to the best of our knowledge.
For MATLAB users, we highlight \texttt{VEMComp} \cite{Frittelli2024}, \texttt{mVEM} \cite{mVEM}, and \texttt{VEMLAB} \cite{Vemlab}.
In the C++ domain, we mention the libraries \texttt{Veamy} \cite{Veamy}, and \texttt{VEM++}.
Finally, the recently developed \texttt{DUNE-VEM} module \cite{DuneVEM} is worth noting, as it is available for both C++ and Python.

The structure of the paper is as follows. After briefly introducing the main notations in Section \ref{sec:notation}, we present the VEM framework in Section \ref{sec:vem} by focusing on the local spaces supported by the library. Section \ref{sec:code} discusses the implementation details, outlining the design of local and global discrete matrices as well as integration and solving techniques available in the library. Finally, Section \ref{sec:examples} reviews selected numerical experiments performed by the authors, where PolyDiM is successfully applied, and presents novel examples showcasing the recently integrated features.

\begin{table}[h]
    \caption{List of the main software available with a framework similar to PolyDiM.}
    \label{tab:software}
    \centering
\begin{tabular}{@{}ccccl@{}}
\toprule
\textbf{Software} & \textbf{Language} & \textbf{Availability} & \textbf{Problem} & \textbf{Website}                       \\ \midrule
\texttt{PolyDiM}           & C++            & Open-Source                  & 2D - 3D          & \href{https://polydim.it/}{https://polydim.it/} \\
\texttt{PolyDiM}           & MATLAB / Python           & Upon Request                 & 2D - 3D          & \href{https://polydim.it/}{https://polydim.it/} \\ \midrule
\texttt{VEMComp}           & MATLAB            & Open-Source                  & 2D - 3D          & \href{https://github.com/massimofrittelli/VEMcomp}{https://github.com/massimofrittelli/VEMcomp} \\ % ref  https://doi.org/10.1007/s11075-024-01919-4
\texttt{mVEM}             & MATLAB               & Open-Source                  & 2D - 3D               & \href{https://github.com/Terenceyuyue/mVEM}{https://github.com/Terenceyuyue/mVEM} \\ % ref https://arxiv.org/abs/2204.01339
\texttt{VEMLAB}             & MATLAB               & Open-Source                  & 2D               & \href{https://camlab.cl/vemlab/}{https://camlab.cl/vemlab/} \\ % ref https://camlab.cl/vemlab/
\texttt{Veamy}             & C++               & Open-Source                  & 2D               & \href{https://camlab.cl/veamy}{https://camlab.cl/veamy} \\ % ref https://doi.org/10.1007/s11075-018-00651-0
\texttt{VEM++}             & C++               & Upon Request                  & 2D - 3D               & \href{https://sites.google.com/view/vembic/home}{https://sites.google.com/view/vembic/home} \\ % ref https://doi.org/10.1007/s11075-025-02059-z
\texttt{DUNE-VEM}          & C++ / Python      & Open-Source                  & 2D - 3D          & \href{https://dune-project.org/sphinx/content/sphinx/dune-fem/vemdemo_descr.html}{https://dune-project.org}         \\ \bottomrule % ref https://doi.org/10.1137/23M1573653
\end{tabular}
\end{table}

\section{Notations}\label{sec:notation}

Throughout this paper, we adopt the following notations. Let $d \in \{1,2,3\}$ be the geometric dimension of the problem.
We denote by $\xx$ and $\yy$ some generic points in $\R^d$ and by $x_i$, with $i \in \{ 1,\dots,d \}$, the $i$-the component of the generic point $\xx$ with respect to Cartesian axes. Let us denote by $E \subset \R^d$ a generic polytope, i.e. a polygon if $d=2$ or a polyhedron if $d=3$. We denote the generic edge of a polytope $E$ by $e$ and, if $d=3$, we use the symbol $F$ to refer to a generic face of the polyhedron $E$. We further denote by $\Nv[E]$, $\Ne[E]$, and $\Nf[E]$ the number of vertices, edges, and faces of $E$, respectively. Moreover, $\vert E \vert$, $\xx_E$ and $h_E = \max_{\xx, \yy \in E} \norm{\xx-\yy}$ represent the measure, the centroid and the diameter of the polygon $E$, respectively. 

Let $\Omega \subset \R^d$ be an open bounded convex polytopal domain and let $\Th$ be a decomposition of $\Omega$ into star-shaped polytopes $E$ that satisfies standard mesh assumptions \cite{BEIRAODAVEIGA20171110}, where we fix, as usual, $h = \max_{E \in \Th} h_E$. 

We define $\Eh[E]$ and $\Fh[E]$ as the set of edges and faces of $E \in \Th$, respectively. Consequently, the set $\Eh[F]$ incorporates all the edges of the face $F$. We use $\partial$ operator to denote the boundary of a polygon or polyhedron. More specifically, if $E$ is a polyhedron, $\partial E$ is the set of faces contained in $\Fh[E]$, while $\partial E$ is the union of edges belonging to $\Eh[E]$ if $E$ is a polygon. Moreover, we set $\Eh = \bigcup_{E\in\Th} \Eh[E]$ and $\Fh = \bigcup_{E\in\Th} \Fh[E]$. We adopt the symbols $\nn_{\partial E}$, $\nn_{F}$ and $\nn_e$ to denote the outward unit normal vectors to the boundary of $E$, to the face $F$ and the edge $e$, respectively. 

We follow standard notations for Sobolev and Hilbert spaces. Let us consider a generic open subset $\genericset \subset \R^d$. Given two scalar functions $p,q \in \leb{2}{\genericset}$, two vector fields $\uu,\ \vv \in \vleb{2}{\genericset}{d}$ and two tensor fields $\bm{T},\ \bm{\sigma} \in \vleb{2}{\genericset}{d \times d}$, we denote by
\begin{align*}
\scal[\genericset]{p}{q} = \int_{\genericset} p q \ ~d\genericset,\quad\scal[\genericset]{\uu}{\vv} = \int_{\genericset} \uu\cdot \vv\ ~d\genericset,\quad \scal[\genericset]{\bm{T}}{\bm{\sigma}} = \int_{\genericset} \bm{T} : \bm{\sigma} \ ~d\genericset, 
\end{align*}
where $\bm{T} : \bm{\sigma} := \sum_{i,j = 1}^n \bm{T}_{ij} \bm{\sigma}_{ij}$. Furthermore, given a generic Sobolev space $\mathcal{H}$, we use the symbol $\norm[\mathcal{H}]{}$ to indicate the norm in $\mathcal{H}$.

Let $\nn_{\Gamma}$ be the outward unit normal vector to the boundary {$\Gamma$} of $\Omega$, we define the following spaces
\begin{gather*}
    \sob[0]{1}{\Omega} = \{v \in \sob{1}{\Omega}: v_{|\Gamma} = 0 \text{ on } \Gamma\},\\
    \Hdiv{\Omega} = \big\{ \vv \in \vleb{d}{\Omega}{d}: \nabla \cdot \vv \in \leb{d}{\Omega}\big\},\\
    \Hdiv[0]{\Omega} = \big\{ \vv \in \Hdiv{\Omega} : \vv \cdot \nn_{\Gamma} = 0 \text{ on } \Gamma \big\}.
\end{gather*}
Additionally, for $d=2$, we define the functional space
\begin{equation*}
    \Hrot{\Omega} = \{\vv\in \vleb{2}{\Omega}{2}: \rot \vv \in \leb{2}{\Omega}\},
\end{equation*}
whereas, for $d=3$, we introduce
\begin{equation*}
    \Hcurl{\Omega} = \Big\{\vv\in \vleb{2}{\Omega}{3}: \curl \vv \in \vleb{2}{\Omega}{3}\Big\}.
\end{equation*}
% In particular, the symbol $\sob{-\frac{1}{2}}{\Gamma}$ denotes the dual space of the Sobolev space $\sob{\frac{1}{2}}{\Gamma}$, whereas the symbol $\dual[\Gamma]{}{}$ denotes the duality pairing between $\sob{-\frac{1}{2}}{\Gamma}$ and $\sob{\frac{1}{2}}{\Gamma}$. 

\subsection{Polynomial spaces}

Let us denote by $\Poly[d]{k}{E}$ the set of the $d$-dimensional polynomials defined on $E$ of degree less or equal to $k \geq 0$ and by $n^d_k = \dim \Poly[d]{k}{E} = \frac{(k+1)\dots(k+d)}{d!}$. For the ease of notation, we further set $\Poly[d]{-1}{E} = \{0\}$ and $n^d_{-1} = 0$ and we use the two natural functions $\ell_{2}: \mathbb{N}^2 \leftrightarrow \mathbb{N}$ and $\ell_{3}: \mathbb{N}^3 \leftrightarrow \mathbb{N}$ that associate:
\begin{equation*}
\begin{aligned}
     \left(  0, 0 \right) \leftrightarrow 1, \quad  \left(  1, 0 \right)                      & \leftrightarrow 2, \quad  \left(  0, 1 \right) \leftrightarrow 3, \quad \left(  2, 0 \right) \leftrightarrow 4, \dots                       \\ 
    \left(  0, 0 ,0 \right) \leftrightarrow 1, \ \left(  1, 0 , 0 \right) \leftrightarrow 2, & \  \left(  0, 1 ,0\right) \leftrightarrow 3, \ \left(  0, 0,1 \right) \leftrightarrow 4,  \ \left(  2, 0,0 \right) \leftrightarrow 5, \dots
    \end{aligned}
    \label{eq:polynomials:ell_monomial} 
\end{equation*}
On each polytope $E$, we define the set of scaled monomials of degree less or equal to $k \geq 0$ as the set
\begin{equation}
    \M[d]{k}{E} = \left\{ \mpoly^{k,d}_{\alpha} = \left(\frac{\xx-\xx_E}{h_E}\right)^{\bm{\alpha}}: \bm{\alpha} = \ell_d(\alpha) \in \mathbb{N}^d, \alpha = 1,\dots, n^d_k\right\},
    \label{eq:polynomials:setscaledmonomial}
\end{equation}
which is a polynomial basis for $\Poly[d]{k}{E}$. Other polynomial bases can be used as well. 
An example is offered by the so-called \textit{MGS basis} described in \cite{Mascotto2018, DassiMascotto2018} and tested by the authors in \cite{Teora2024}. Therefore, for the sake of generality, we introduce the generic set of polynomial basis functions $\mathcal{P}_{k}^d(E) = \{\ppoly^{k,d}_{\alpha}\}_{\alpha=1}^{n^d_k}$ for the polynomial space $\Poly[d]{k}{E}$ defined over a generic polytope $E$. 

Starting from a generic polynomial basis for $\Poly[d]{k}{E}$, an easily computable basis $\{\pp_I^{k,d}\}_{I=1}^{dn_k^d}$ for $\VPoly[d]{k}{E}{d}$ can be built as
\begin{equation}
    \pp_I^{k,d} =\begin{cases}
    \begin{bmatrix} \ppoly^{k,d}_I\\ \vdots \\ 0\end{bmatrix} & I=1,\dots,n^d_k,\\ 
    \cdots \\
    \begin{bmatrix} 0 \\ \vdots \\ \ppoly^{k,d}_{I-(d-1)n_k}\end{bmatrix} & I=(d-1)n^d_k+1,\dots,dn_k.
    \end{cases}
    \label{eq:polynomials:vector_monomials}
\end{equation}

Furthermore, we introduce the vector-polynomial space
\begin{equation}
    \GPoly{k}{\nabla}{E} = \nabla \Poly[d]{k+1}{E} := \myspan \big\{\g^{\nabla,k}_{\alpha}\big\}_{\alpha = 1}^{n_k^{\nabla}}
    \label{eq:polynomials:GNablak}
\end{equation}
and its orthogonal complement $\GPoly{k}{\bot}{E} = \myspan \big\{\g^{\bot,k}_{\alpha}\big\}_{\alpha = 1}^{n_k^{\bot}}$ in $\VPoly[d]{k}{E}{d}$, which satisfies
\begin{equation}
    \VPoly[d]{k}{E}{d} = \GPoly{k}{\nabla}{E} \bigoplus \GPoly{k}{\bot}{E},
    \label{eq:polynomials:complemCondition}
\end{equation}
where $\bigoplus$ is the direct sum operator, and
\begin{equation*}
    \dim \VPoly[d]{k}{E}{d} =d n_k^d,
\end{equation*}
\begin{equation*}
   n^{\nabla}_k = \dim {\GPoly{k}{\nabla}{E}} = n_{k+1}^d - 1,
\end{equation*}
\begin{equation*}
   n^{\bot}_{k} = \dim \GPoly{k}{\bot}{E} = d n_k^d - n^{\nabla}_k.
\end{equation*}

Basis functions for the $\GPoly{k}{\nabla}{E}$ space can be written as
\begin{equation*}
    \g^{\nabla,k}_{\alpha} = \nabla \ppoly^{k+1,d}_{\alpha + 1} = \sum_{I = 1}^{d n_k} \TT{\nabla,k}_{\alpha I} \pp_I^{k,d} \quad \forall \alpha = 1,\dots, n_k^{\nabla},
\end{equation*}
where $\TT{\nabla,k} \in \R^{n^{\nabla}_k \times dn_k}$ is the coefficient matrix that expresses gradients of the polynomial functions $ \{\ppoly^{k+1,d}_{\alpha}\}_{\alpha=2}^{n^d_{k+1}}$ with respect to the polynomial basis $\{\pp^{k,d}_{I}\}_{I=1}^{dn_k}$ of $\VPoly[d]{k}{E}{d}$.
Moreover, each $\g^{\bot,k}_{\alpha}$ function can be written as
\begin{equation*}
    \g^{\bot,k}_{\alpha}  = \sum_{I = 1}^{d n_k} \TT{\bot,k}_{\alpha I} \pp_I^k \quad \forall \alpha = 1,\dots, n_k^{\bot},
\end{equation*}
where $\TT{\bot,k} \in \R^{n^{\bot}_k \times dn_k}$ is the matrix whose rows define an \textit{Euclidean} orthonormal basis for the nullspace of the $\TT{\nabla,k}$ matrix. Thus, by considering the singular value decomposition of $\TT{\nabla,k} = \USVD\ \SSVD \left(\VSVD\right)^T$, we can define $\TT{\bot,k}$ as 
\begin{equation*}
    \TT{\bot,k} = \left[\VSVD(:,n^{\nabla}_k+1:dn^d_k)\right]^T,
\end{equation*}
where $\VSVD(:,n^{\nabla}_k+1:dn^d_k)$ is the sub-matrix of $\VSVD$ made up of all its rows and of the columns running from the $(n_{k}^{\nabla}+1)$-th to the $dn^d_{k}$-th. As a consequence,
\begin{equation*}
    \TT{\nabla,k} \left(\TT{\bot,k}\right)^T = \Omatrix.
\end{equation*}

\begin{remark}
    We observe that the space $\GPoly{k}{\perp}{E}$ could be replace by each complement $\GPoly{k}{\bigoplus}{E}$ of $\GPoly{k}{\nabla}{E}$ in $\VPoly[d]{k}{E}{d}$. For example, $\GPoly{k}{\bigoplus}{E} = \xx^{\perp} \Poly{k-1}{E}$ with $\xx^{\perp} := [x_2, -x_1]^T$ if $d=2$, whereas it can be written as $\GPoly{k}{\bigoplus}{E} = \xx \wedge \VPoly[3]{k-1}{E}{3}$, if $d=3$ \cite{Dassi2019}.
\end{remark}

To perform efficient matrix-based computations in PolyDiM, the evaluations of such polynomial basis functions at given evaluation points are stored in a Vandermonde matrix, whose $(i,j)$-th entry represents the evaluation of the $j$-th basis function at the $i$-th evaluation point. 
These matrices are built through functions available in the namespace \lstinline[language=C++, style=mystyle]{Polydim::VEM::Utilities}. 

\section{The Virtual Element Framework}\label{sec:vem}

In this section, we introduce the Virtual Element Method framework, focusing on the definition and structure of the local spaces implemented in the PolyDiM library. We also describe the main VEM projection operators and the stabilization terms, highlighting their role in ensuring the method consistency and stability properties.

In the following, we will refer to the \textit{primal formulation} when the discretization involves only a single variable $u$ (or $\uu$ when considering vector-problems). On the other hand, the \textit{mixed formulation} will refer to a discretization involving two variables $(\uu,p)$, to which we refer as the velocity field and the pressure variable, respectively. Finally, \textit{divergence-free formulation} refers to the discretization of a problem involving two variables $(\uu,p)$, again called the velocity and the pressure fields, whose related velocity space is point-wise divergence-free \cite{DaVeigaLovadina2017}.

For any $k \in \N$ and for each polytope or polytopal face $\genericelement$, we introduce the following polynomial projectors:
\begin{itemize}
    \item the $L^2$-projection $\proj{0,\genericelement}{k} : \leb{2}{\genericelement} \to \Poly[d]{k}{\genericelement}$, which is defined for all $v \in \leb{2}{\genericelement}$ by
    \begin{equation}
        \scal[\genericelement]{v - \proj{0, \genericelement}{k} v}{\ppoly} = 0,\quad \forall \ppoly \in \Poly[d]{k}{\genericelement},
        \label{eq:proj:Pi0k}
    \end{equation} 
    with natural extension $\proj{0,\genericelement}{k} : \vleb{2}{\genericelement}{d} \to \VPoly[d]{k}{\genericelement}{d}$ for vector fields $\vv \in \vleb{2}{\genericelement}{d}$.
    \item the $H^1$-seminorm projection $\proj{\nabla,\genericelement}{k} : \sob{1}{\genericelement} \to \Poly[d]{k}{\genericelement}$, such that for any $v \in \sob{1}{\genericelement}$
    \begin{equation}
    \resizebox{0.8\hsize}{!}{
        $\begin{cases}
            \scal[\genericelement]{\nabla v - \nabla \proj{\nabla, \genericelement}{k} v }{\nabla \ppoly} = 0 & \forall \ppoly \in \Poly[d]{k}{\genericelement},\\
            P_0(\proj{\nabla,E}{k}v - v ) = 0, 
        \end{cases} \quad \text{where} \quad P_0(v) = \begin{cases}
        \int_{\partial \genericelement} v & \text{if } k =1,\\
        \int_{\genericelement}  v & \text{if } k \geq 2.
        \end{cases}$}
        \label{eq:proj:PiNabla}
    \end{equation}
    with clear extension $\proj{\nabla,\genericelement}{k} : \vsob{1}{\genericelement}{d} \to \VPoly[d]{k}{\genericelement}{d}$ for vector fields $\vv \in \vsob{1}{\genericelement}{d}$.
\end{itemize}

\subsection{Primal $\sob{1}{\Omega}$-Conforming}\label{sec:PCC}

Let us consider the following Poisson problem 
\begin{equation}
    \begin{cases}
        \nabla\cdot \left(-\D \nabla u\right) = f & \text{in } \Omega,\\
        u = 0  & \text{on } \Gamma,
    \end{cases}
    \label{eq:primal:continuousmodelproblem}
\end{equation}
where $\D \in \R^{d \times d}$ is a symmetric positive-definite tensor and $f$ is a sufficiently regular scalar field.
Here, for the sake of simplicity, we impose homogeneous Dirichlet boundary conditions, but general boundary conditions as well as more general elliptic problems can be handled with PolyDiM.

Setting $\VP := \sob[0]{1}{\Omega}$, the variational formulation of Problem \eqref{eq:primal:continuousmodelproblem} reads as: \textit{Find} $u \in \VP$, \textit{such that}
\begin{equation}
    \dbilin{u}{v} =  \scal[\Omega]{f}{q} \qquad \forall v \in  \VP
    \label{eq:primal:variationalfomrulation}
\end{equation}
where 
\begin{equation}
    \dbilin{u}{v} := \scal[\Omega]{\D \nabla u}{ \nabla v} \quad \forall u,v \in \VP.
    \label{eq:primal:continuous_dbilin}
\end{equation}

\subsubsection{The Local Space and the \lstinline[language=C++, style=mystyle]{Polydim::VEM::PCC} classes}

Let $k \geq 1$, for each polygon or polygonal face $\genericelement$, we introduce the space
\begin{equation*}
    \Bk{k}{\partial \genericelement} = \left\{v \in \con{0}{\partial \genericelement}: v_{|e}\!\in\!\Poly[1]{k}{e} \ \forall e \in \Eh[\genericelement] \right\},
\end{equation*}
which allows us to define the two-dimensional local primal virtual element space on the element $E$ as
\begin{equation}
\begin{aligned}
    \VPh[E]{k,\ell}{2} = \big\{v \in \sob{1}{E}:&(i)\ \ v_{|\partial E}\!\in\!\Bk{k}{\partial E},\\ 
    &(ii)\ \Delta v \in \Poly[2]{k+\ell}{E},\\  
    &(iii) \scal[E]{v - \proj{\nabla,E}{k} v}{\ppoly} = 0 \ \forall \ppoly \in \Poly[2]{k + \ell}{E} \setminus\Poly[2]{k-2}{E}  \big\}.
\end{aligned}
\label{eq:primal:twolocalvirtualspace}
\end{equation}
The set $\Poly[d]{k+\ell}{E} \setminus\Poly[d]{k-2}{F}$ denotes the union of set of $d$-dimensional homogeneous polynomials defined on $E$ from degree $k-1$ to $k + \ell$. If $d=3$, for each polyhedron $E \in \Th$, we further consider the boundary space
\begin{equation}
    \mathbb{U}_{k,\ell}(\partial E) = \left\{v \in \con{0}{\partial E}: v_{|F}\!\in\!\VPh[F]{k,\ell}{2} \ \forall F \in \Fh[E]  \right\},
    \label{eq:primal:virtualfacespace}
\end{equation} 
and then we define three-dimensional local primal virtual element space as
\begin{equation}
\begin{aligned}
\VPh[E]{k,\ell}{3} = \big\{v \in \sob{1}{E}: &(i)\ \ v_{| \partial E} \!\in\! \mathbb{U}_{k,\ell}(\partial E),\\
&(ii)\ \Delta v \in \Poly[3]{k + \ell}{E},\\
&(iii) \scal[E]{v - \proj{\nabla,E}{k} v}{\ppoly} = 0 \ \forall \ppoly \in \Poly[3]{k+\ell}{E} \setminus\Poly[3]{k-2}{E}\}.
\end{aligned}
\label{eq:primal:threelocalvirtualspace}
\end{equation}
We want to remark that the Property $(iii)$ which defines functions in the local virtual element spaces is usually called the \textit{enhancement property}.

\begin{figure}[!h]
    \centering
    \begin{subfigure}{0.33\linewidth}
        \includegraphics[width=0.9\linewidth]{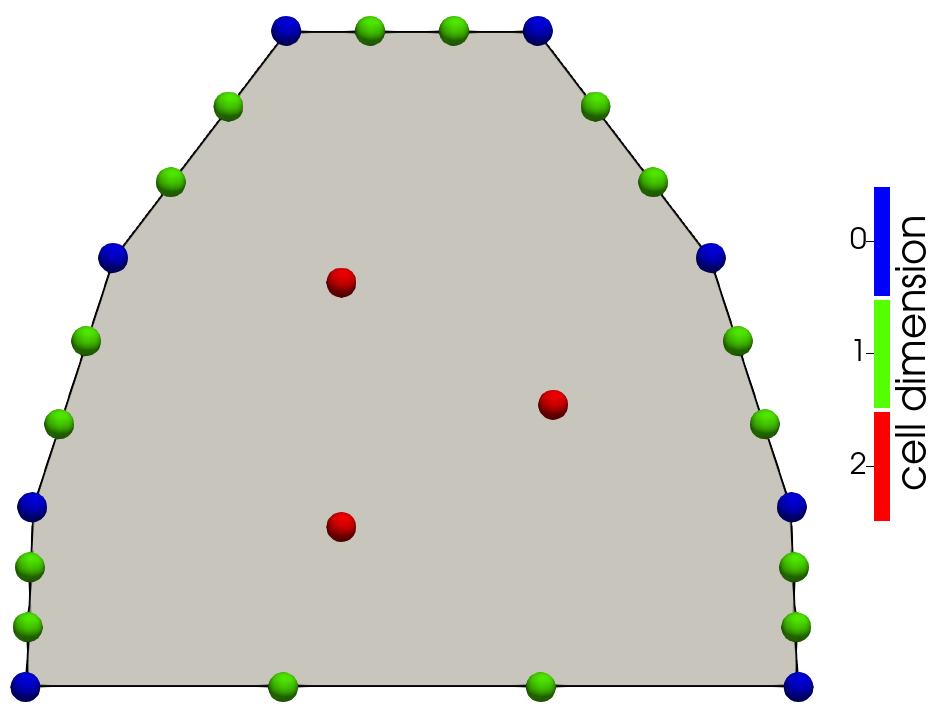}
        \caption{}
        \label{fig:primal:dofs}
    \end{subfigure}
    \begin{subfigure}{0.33\linewidth}
        \includegraphics[width=0.9\linewidth]{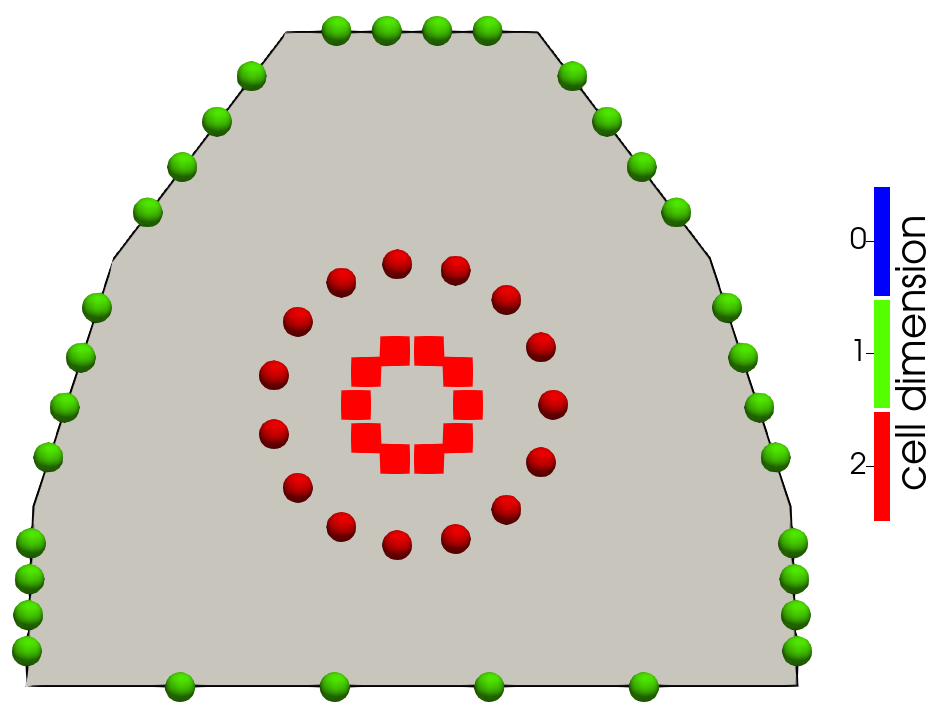}
        \caption{}
        \label{fig:mixed:dofs}
    \end{subfigure}
    \begin{subfigure}{0.33\linewidth}
        \includegraphics[width=0.9\linewidth]{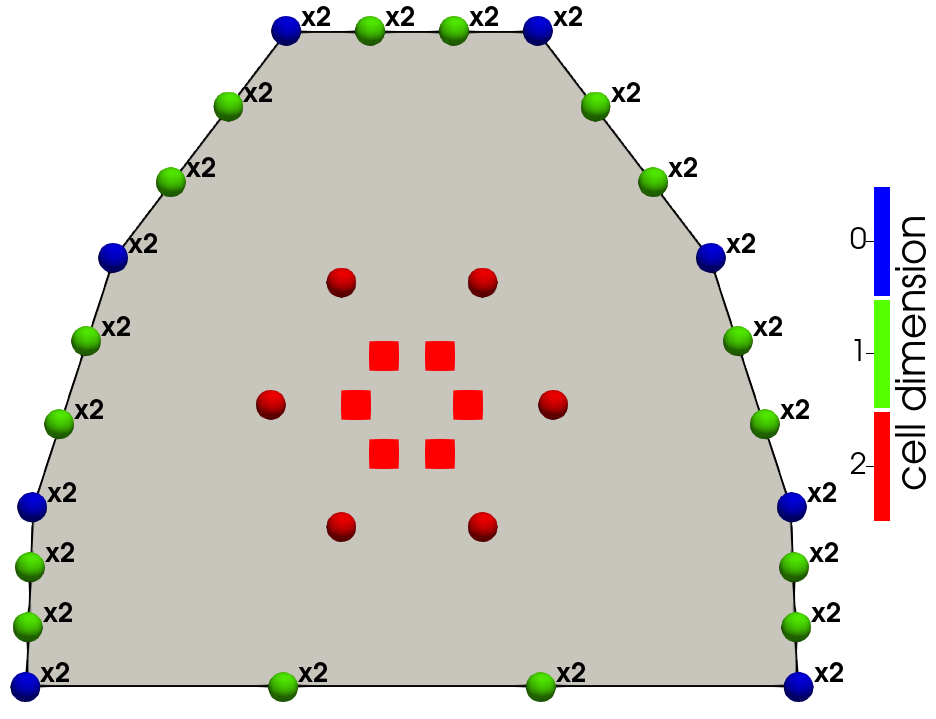}
        \caption{}
        \label{fig:divfree:dofs}
    \end{subfigure}
    \caption{Degrees of freedom for $k=3$ and $d=2$. Left: DOFs associated with the local primal space $\VPh{k}{d}$ defined in \eqref{eq:primal:twolocalvirtualspace}. Center: Degrees of freedom for the local mixed spaces. Circles refer to the velocity space $\VMh[E]{k}{d}$, defined in \eqref{eq:mixed:VemVelocitySpace_2D}, whereas squares refer to the pressure space $\Qh[E]{k}{d}$. Right: Degrees of freedom for the local divergence-free spaces. Circles refer to the velocity space $\VMh[E]{k}{d}$, defined in \eqref{eq:divfree:discrete_local_velocity_space_stokes}, whereas squares refer to the pressure space $\Qh[E]{k}{d}$ defined in \eqref{eq:divfree:discrete_local_pressure_space_stokes}. The symbol $\times 2$ near a marker denotes the multiplicity of the related degree of freedom.}
    \label{fig:dofs}
\end{figure}

For each $v \in \VPh[E]{k,\ell}{d}$, a possible choice of the local Degrees of Freedom (DOFs) are: 
\begin{itemize}
\item \textbf{Vertex DOFs}: the values of $v$ at each vertex of $E$.
\item \textbf{Edge DOFs}: the values of $v$ at the $k-1$ internal Gauss-Lobatto quadrature points on each edge $e \in \Eh[E]$.
\item \textbf{Face DOFs}: if $d=3$, the internal moments on each face $F\in\Fh[E]$ defined as
\begin{equation}
{\frac{1}{\vert F \vert}}\int_F v \ \ppoly \qquad \forall \ppoly \in \mathcal{P}_{k-2}^2(F).
\label{eq:primal:face_DOF}
\end{equation}
\item \textbf{Internal DOFs}: the internal moments:
\begin{equation}
 {\frac{1}{\vert E \vert}}\int_E v \ \ppoly \qquad \forall \ppoly \in \mathcal{P}_{k-2}^d(E).
\label{eq:primal:internal_DOF}
\end{equation}
\end{itemize}
The total number of local degrees of freedom is
\begin{equation}
    \Ndof[E] := \begin{cases}
    \Nv[E] + \Ne[E] (k-1) + \frac{(k-1)k}{2} & \text{if } d=2,\\
    \Nv[E] + \Ne[E] (k-1) + \Nf[E] \frac{(k-1)k}{2} + \frac{(k-1)k(k+1)}{6} & \text{if } d=3.
    \end{cases}
\end{equation}
A graphical illustration of such degrees of freedom for $k=3$ in two-dimensions is shown in Figure \ref{fig:primal:dofs}. It is important to note that the DOF positions shown are only illustrative. Indeed, many DOFs are not associated with specific points on the element. 
\begin{proposition}[\cite{LBe16, LBe13, Ahmad2013}]\label{prop:primal:projection}
    Let $\VPh[E]{k,\ell}{d}$ be the space defined in \eqref{eq:primal:twolocalvirtualspace} for $d=2$ and in \eqref{eq:primal:threelocalvirtualspace} for $d=3$. The local DOFs allow us to compute exactly (up to machine precision) the element projections
    \begin{equation*}
        \proj{0,E}{k-2} : \VPh[E]{k,\ell}{d} \to \Poly[d]{k-2}{E},\quad \proj{\nabla,E}{k} : \VPh[E]{k,\ell}{d} \to \Poly[d]{k}{E}, \quad \vproj{0,E}{k-1} : \nabla \VPh[E]{k,\ell}{d} \to \VPoly[d]{k-1}{E}{d}
    \end{equation*}
     as defined in \eqref{eq:proj:Pi0k} and \eqref{eq:proj:PiNabla}, i.e. we are able to compute the polynomials $\proj{0,E}{k-2} v_h$, $\vproj{\nabla,E}{k} v_h$ and $\vproj{0}{k-1} \nabla v_h$, for any $v_h \in \VPh[E]{k}{d}$, using only the information given by the local degrees of freedom of $v_h$. Moreover, by exploiting the enhancement property, we are also able to compute the higher-order projections
    \begin{equation*}
        \proj{0,E}{\ast} : \VPh[E]{k,\ell}{d} \to \Poly[d]{\ast}{E} \quad \forall \ast = k-1,\dots,k+\ell.
    \end{equation*}
\end{proposition}

PolyDiM offers the possibility to choose among different sets of degrees of freedom and different approaches to build the polynomial projections of virtual element functions. In particular, on the current release the reader can find:
\begin{itemize}
    \item The monomial approach, that exploits the monomial basis defined in \eqref{eq:polynomials:setscaledmonomial} in the definition of local degrees of freedom \eqref{eq:primal:face_DOF} and \eqref{eq:primal:internal_DOF} as well as in the definition of local projections. The corresponding approach is implemented following \cite{DaVeiga2016} and the local space can be built exploiting classes \lstinline[language=C++, style=mystyle]{Polydim::VEM::PCC::VEM_PCC_2D_LocalSpace} if $d=2$ and \lstinline[language=C++, style=mystyle]{Polydim::VEM::PCC::VEM_PCC_3D_LocalSpace} if $d=3$.
    \item The orthonormal approach introduced in \cite{Mascotto2018, DassiMascotto2018} and tested by the authors in \cite{Teora2024}. This approach guarantees well-conditioned local and global matrices for the higher-values of the local polynomial degree $k$. The related local space can be built exploiting \lstinline[language=C++, style=mystyle]{Polydim::VEM::PCC::VEM_PCC_2D_Ortho_LocalSpace} if $d=2$ and \lstinline[language=C++, style=mystyle]{Polydim::VEM::PCC::VEM_PCC_3D_Ortho_LocalSpace} if $d=3$.
    \item The novel (B-F) inertial approach introduced by the authors in \cite{Teora2024} that reveals to be robust in the presence of badly-shaped polytopes. The corresponding local spaces can be built exploiting the class \lstinline[language=C++, style=mystyle]{Polydim::VEM::PCC::VEM_PCC_2D_Inertia_LocalSpace} if $d=2$ and \lstinline[language=C++, style=mystyle]{Polydim::VEM::PCC::VEM_PCC_3D_Inertia_LocalSpace} if $d=3$.
\end{itemize}

Let us introduce the set $\{\varphi_{j,E}\}_{j=1}^{\Ndof[E]}$ of Lagrange basis functions related to the aforementioned local DOFs, where the numbering $j$ counts in order the vertex, the edge, the face, and the internal DOFs. 

In the virtual element framework, basis functions are not polynomial as for the FEM and are not known in a closed-form in the interior of an element $E \in \Th$. Thus, to solve a partial differential equation using VEM, we must resort to the projection operators defined in Proposition \ref{prop:primal:projection}. All the information needed to compute such projections is stored in the data structure \lstinline[language=C++, style=mystyle]{Polydim::VEM::PCC::VEM_PCC_dD_LocalSpace_Data}.

The above-mentioned classes offer user-friendly interfaces that allows us to retrieve the values of the projection of the local basis functions $\{\varphi_{j,E}\}_{j=1}^{\Ndof[E]}$ and their derivatives in a given set of points. The list of available projection operators is defined in the enumeration \lstinline[language=C++, style=mystyle]{Polydim::VEM::PCC::ProjectionTypes} and coincides with projectors defined in Proposition \ref{prop:primal:projection}.

We highlight that, since virtual element functions are known polynomial on $\partial E$, we use the exact expression of the Lagrange basis functions to evaluate boundary contributions.

\subsubsection{The Discrete Primal Variational Formulation}\label{sec:PCC:discreteformulation}

We define the global Virtual Element space as \cite{BEIRAODAVEIGA20171110}:
\begin{equation}
\VPh{k}{} = \{v \in \VP \cap \con{0}{\overline{\Omega}}: v \in \VPh[E]{k, 0}{d} \ \forall E \in \Th~ \big\}.
\label{eq:primal:viartualspace}
\end{equation}

The continuous bilinear form \eqref{eq:primal:continuous_dbilin} can be split according to the tessellation $\Th$ as
\begin{equation*}
\dbilin{u}{v} = \sum_{E \in \Th} \dbilin[E]{u}{v},\quad \dbilin[E]{u}{v} = \scal[E]{\D \nabla u}{ \nabla v}\quad \forall u,v \in \VP. 
\end{equation*}
In general, we are not able to compute the quantity
\begin{equation*}
\dbilin[E]{u_h}{v_h} = \scal[E]{\D \nabla u_h}{ \nabla v_h}\quad \forall u_h,\ v_h \in \VPh[E]{k}{d}, 
\end{equation*}
since we do not know the virtual element functions in a closed form in the interior of each element $E\in \Th$.
To overcome this issue, the main idea of the Virtual Element Method is to substitute the continuous bilinear form with a \textit{computable} discrete counterpart $\dbilinh[E]{}{}: \VPh[E]{k}{d} \times  \VPh[E]{k}{d} \to \R$ which satisfies the two following properties \cite{LBe13}:
\begin{itemize}
\item \textit{Consistency}: For all $p \in \Poly[d]{k}{E}$ and for all $v_h \in \VPh[E]{k}{d}$
\begin{equation*}
\dbilinh[E]{p}{v_h} = \dbilin[E]{p}{v_h}.
\end{equation*}
\item \textit{Stability}: There exist two positive constants $\alpha_{\ast},\ \alpha^{\ast}$ independent of $h$ such that
\begin{equation}
    \alpha_{\ast} \dbilin[E]{v}{v} \leq \dbilinh[E]{v}{v} \leq \alpha^{\ast} \dbilin[E]{v}{v},\quad \forall v \in \VPh[E]{k}{d}:\ \proj{\nabla,E}{k}v = 0.
    \label{eq:primal:stab_property}
\end{equation}
\end{itemize}
To build a discrete bilinear form that satisfies the consistency and stability properties, the local continuous bilinear form is first split as 
\begin{equation}
\dbilin[E]{u_h}{v_h} = \dbilin[E]{\proj{\nabla,E}{k}u_h}{\proj{\nabla,E}{k}v_h} + \dbilin[E]{(I-\proj{\nabla,E}{k})u_h}{(I-\proj{\nabla,E}{k})v_h},
\label{eq:primal:cont_vem_bilinear_form}
\end{equation}
where the equality is due to the orthogonality of $\proj{\nabla,E}{k}$ with respect to the scalar product induced by $\dbilin[E]{}{}$. The first term in the right-hand side of \eqref{eq:primal:cont_vem_bilinear_form} is computable thanks to the definition of the local degrees of freedom, whereas the second one could be approximated by any \textit{computable} symmetric positive definite bilinear form $\stab[E]{}{}$ that satisfies the stability property \eqref{eq:primal:stab_property}.
In \cite{LBe16}, the local bilinear form $\dbilin[E]{\proj{\nabla,E}{k}u_h}{\proj{\nabla,E}{k}v_h}$ is replaced by 
\begin{equation*}
    \scal[E]{\D\proj{0,E}{k-1}\nabla u_h}{\proj{0,E}{k-1}\nabla v_h}
\end{equation*}
to avoid loss of accuracy for the high orders of the method when dealing with variable coefficients. We observe that, when $k=1$, the two choices coincide with each other.
Finally, the local virtual discrete bilinear form is 
\begin{equation*}
\dbilinh[E]{u_h}{v_h} := \scal[E]{\D\proj{0,E}{k-1}\nabla u_h}{\proj{0,E}{k-1}\nabla v_h} + \stab[E]{(I-\proj{\nabla,E}{k})u_h}{(I-\proj{\nabla,E}{k})v_h},
\end{equation*}
which is computable and satisfies the consistency and the stability property \cite{LBe13}.
Now, let us define $\dof_i^E$, for each $i=1,\dots,\Ndof[E]$ and each $E\in \Th$, as the operator that associates with each sufficiently smooth function $\varphi$ its $i$-th local degree of freedom $\dof_i^E(\varphi)$. A standard choice for the stabilization term is given by the \textit{dofi-dofi} stabilization term 
\begin{equation}
    \stab[E]{u_h}{v_h} = h^{d-2} \sum_{i=1}^{\Ndof[E]} \dof_i^E(u_h) \dof_i^E(v_h).
    \label{eq:primal:dofidofi}
\end{equation}
We observe that, when dealing with more general elliptic equations, this stabilization is usually pre-multiplied by a constant $C_{s}$, which accounts for the magnitude of the diffusion coefficients.
Other stabilization methods have been proposed in the literature, which may take integral forms \cite{LBe17} or be a variant of the dofi-dofi stabilization, such as the $D$-recipe version introduced in \cite{BEIRAODAVEIGA20171110}. In particular, the $D$-recipe form aims to prevent the stabilization from becoming too small in magnitude with respect to the consistency term when high-order methods are considered. 

Finally, the virtual element discretization of problem \eqref{eq:primal:variationalfomrulation} reads as: \textit{Find $u_{h} \in \VPh{k}{}$ such that:} 
\begin{equation}
    \sum_{E \in \Th} \dbilinh[E]{u_h}{v_h}  = \sum_{E \in \Th} \scal[E]{f}{\proj{0,E}{k-1} v_h} \quad \forall v_h \in \VPh{k}{}.
\label{eq:primal:vem_problem_discretization}
\end{equation}
The problem \eqref{eq:primal:vem_problem_discretization} has a unique solution $u_{h} \in \VPh{k}{}$ and, for $h$ sufficiently small, the following a priori error estimates hold true
\begin{equation}
    \norm[\leb{2}{\Omega}]{u - u_h} = O( h^{k+1}),\quad 
    \norm[\vleb{2}{\Omega}{d}]{ \nabla u -  \nabla u_h} = O( h^{k}).
\end{equation}

\subsubsection{A simple extension to deal with vector-problems}\label{sec:elastic}

A simple way to deal with vector-problems, as the elasticity problems, is to build the local space for the primal vector-variable $\uu \in \VM :=\vsob[0,\Gamma_D]{1}{\Omega}{d}$, with $\uu = \begin{bmatrix}
    u_1\\
    \dots\\
    u_d
\end{bmatrix}$, as the cartesian product of the corresponding scalar local space \cite{DaVeigaBrezzi2013}. Specifically,
we define the virtual element local space as $\VMh[E]{k}{d} = [\VPh[E]{k}{d}]^d$.
By proceeding in this way, the number of local degrees of freedom for a vector-valued function $\vv \in \VMh[E]{k}{d}$ becomes $d \Ndof[E]$.

For instance, the virtual element discretization for the linear elastic problem with homogeneous Dirichlet boundary condition reads as: \textit{Find $\uu_h \in \VMh{k}{}$ such that}
\begin{equation}
\begin{aligned}
    &\sum_{E \in \Th} \Big[\scal[E]{2\mu \proj{0,E}{k-1} \bm{\epsilon}(\uu_h)}{\proj{0,E}{k-1}\bm{\epsilon}(\vv_h)} + \mu \stab[E]{(I - \proj{\nabla,E}{k})\uu_h}{(I - \proj{\nabla,E}{k})\vv_h}\\
    &\qquad\qquad +\scal[E]{\lambda \proj{0,E}{k-1} \div \uu_h}{\proj{0,E}{k-1}\div \vv_h} \Big] =  \sum_{E \in \Th} \scal[E]{\ff}{\proj{0,E}{k-1} \vv_h} \quad \forall \vv_h \in \VMh{k}{}
\end{aligned}
\label{eq:primal:vem_elasticity}
\end{equation}
where $\bm{\epsilon}$ denotes the symmetric gradient, $\lambda$ and $\mu$ are positive coefficients (Lamé coefficients), while $\ff$ is a vector-
valued function belonging to $\vleb{2}{\Omega}{d}$. Moreover, $\forall \vv \in \vsob{1}{E}{d}$, we define \cite{DaVeigaBrezzi2023}
\begin{gather*}
    \proj{0,E}{k-1} \bm{\epsilon}(\vv) = \frac{1}{2}\left(\proj{0,E}{k-1} \nabla \vv + (\proj{0,E}{k-1} \nabla \vv)^T\right), \\
    \proj{0,E}{k-1} \div \vv = \proj{0,E}{k-1} \tr (\nabla \vv) = \tr (\proj{0,E}{k-1} \nabla \vv),
\end{gather*}
where $\tr(\cdot)$ denotes the trace of a tensor-field.

\subsection{Mixed $\Hdiv{\Omega}$-conforming}\label{sec:MCC}

To deal with the mixed variational formulation, we define
\begin{equation*}
    p = u, \quad \K = \D^{-1}, \quad \bbeta = \K \bb,
\end{equation*}
and we re-write problem \eqref{eq:primal:continuousmodelproblem} as
\begin{equation}
    \begin{cases}
        \K \uu =  -\nabla p  & \text{ in } \Omega,\\
        \nabla \cdot \uu  = f & \text{ in } \Omega,\\
        p = 0 & \text{ on } \Gamma.
    \end{cases}
    \label{eq:mixed:continuouspmodelroblem}
\end{equation}
Setting
\begin{equation}
    \VM := \Hdiv{\Omega},\quad \Q := \leb{2}{\Omega},
    \label{eq:mixed:continuousspaces}
\end{equation}
the mixed variational formulation of \eqref{eq:primal:continuousmodelproblem} reads as: \textit{Find } $\uu \in \VM, \textit{ and } p \in \Q$ \textit{ such that }
\begin{equation}
    \begin{cases}
        \dbilin{\uu}{\vv} -  \scal[\Omega]{p}{\nabla \cdot \vv} = 0 & \forall \vv \in \VM \\
        \scal[\Omega]{\nabla \cdot \uu}{ q} = \scal[\Omega]{f}{q} & \forall q \in  \Q
    \end{cases}
    \label{eq:mixed:variationalfomrulation}
\end{equation}
where 
\begin{equation}
    \dbilin{\uu}{\vv} := \scal[\Omega]{\K \uu}{\vv}\quad \forall \uu,\vv \in \VM.
    \label{eq:mixed:continuous_dbilin}
\end{equation}

\subsubsection{The Local Space and the \lstinline[language=C++, style=mystyle]{Polydim::VEM::MCC} classes}

Following \cite{secondMixed}, for any integer $k\geq 0$, if $d=2$, we define the local mixed virtual element space for the velocity variable $\uu$ as
\begin{equation}
    \begin{aligned}
        \nonumber \VMh[E]{k}{2} = \big\{ \vv_h \in \Hdiv{E} \cap \Hrot{E} \text{ s.t. } &(i)\ \ \vv_h \cdot \nn_{e} \in \Poly[2]{k}{e}\ \forall e \in \Eh[E], \\
        &(ii)\ \nabla \cdot \vv_h \in \Poly[2]{k}{E},\\
        &(iii)\ \rot \vv_h \in \Poly[2]{k-1}{E} \big\},
    \end{aligned}
    \label{eq:mixed:VemVelocitySpace_2D}
\end{equation}
whereas, for $d=3$, we define the space as \cite{HdivHcurl}
\begin{equation}
    \begin{aligned}
        \nonumber \VMh[E]{k}{3} = \big\{ \vv_h \in \Hdiv{E} \cap \Hcurl{E} \text{ s.t. } &(i)\ \ \vv_h \cdot \nn_{f} \in \Poly[3]{k}{f}\ \forall f \in \Fh[E], \\
        &(ii)\ \nabla \cdot \vv_h \in \Poly[3]{k}{E},\\
        &(iii)\ \curl \vv_h \in \curl \VPoly[3]{k}{E}{3} \big\}.
    \end{aligned}
    \label{eq:mixed:VemVelocitySpace_3D}
\end{equation}
The following set of local degrees of freedom is unisolvent for $\VMh[E]{k}{d}$ (see \cite{DaVeiga2016,HdivHcurl} for more details and Figure \ref{fig:mixed:dofs} for a graphical illustration): for each $\vv_h \in \VMh[E]{k}{d}$, 
\begin{itemize}
    \item \textbf{Boundary DOFs}: if $d=2$, the values of $\vv_h \cdot \nn_e$ at the $k+1$ Gauss quadrature points internal on each edge $e \in \Eh[E]$, 
    or, more generally,
    \begin{equation}
        \begin{cases}
            \frac{1}{\vert e \vert} \int_e \vv_h \cdot \nn_e\ p\ \quad \forall p \in \mathcal{P}_k^1(e)\quad \forall e \in \Eh[E] & \text{if } d=2,\\
            \frac{1}{\vert f \vert} \int_f \vv_h \cdot \nn_f\ p \quad \forall p \in \mathcal{P}_k^2(f)\quad \forall f \in \Fh[E] & \text{if } d=3.
        \end{cases}
        \label{eq:mixed:boundaryDOFs}
    \end{equation}
    We note that this choice automatically ensures the continuity of the flux $\vv_h\cdot \nn$ across two adjacent elements.
    \item \textbf{Internal $\nabla$ DOFs}:
    \begin{equation}
        \frac{1}{\vert E \vert} \int_E \vv_h \cdot \g^{\nabla,k-1}_{\alpha}\quad \forall \alpha = 1,\dots, n_{k-1}^{\nabla}.
        \label{eq:mixed:DOFnabla}
    \end{equation}
    \item \textbf{Internal $\bot$ DOFs}:
        \begin{equation}
        \frac{1}{\vert E \vert} \int_E \vv_h \cdot \g^{\bot,k}_{\alpha}\quad \forall \alpha = 1,\dots, n_{k}^{\bot}.
        \label{eq:mixed:DOFbigoplus}
    \end{equation}
\end{itemize}

The dimension of the local mixed virtual element space is given by
\begin{equation}
 \Ndof[E] = \dim\VMh[E]{k}{d} = N^{\ast}_E n_k^{d-1} +  n_{k-1}^{\nabla} + n_{k}^{\bot},\quad \text{with}\quad \ast = \begin{cases}
    e & \text{if } d=2,\\
    f & \text{if } d =3.
\end{cases}   
\end{equation}

\begin{proposition}[\cite{secondMixed, HdivHcurl}]
    Let $\VMh[E]{k}{d}$ be the space defined in \eqref{eq:mixed:VemVelocitySpace_2D} or in \eqref{eq:mixed:VemVelocitySpace_3D}. The DOFs \eqref{eq:mixed:boundaryDOFs}, \eqref{eq:mixed:DOFnabla} and \eqref{eq:mixed:DOFbigoplus} allow us to compute exactly (up to machine precision) the element projection
    \begin{equation*}
        \vproj{0,E}{k} : \VMh[E]{k}{d} \to \VPoly[d]{k}{E}{d}
    \end{equation*}
     as defined in \eqref{eq:proj:Pi0k}, i.e. we are able to compute the polynomial $\vproj{0,E}{k} \vv_h$, for any $\vv_h \in \VMh[E]{k}{d}$, using only the information given by the local degrees of freedom of $\vv_h$.
\end{proposition}
We further introduce the local mixed virtual element space $\Qh[E]{k}{d}$ for the pressure variable $p$ as the space of polynomials $\Poly[d]{k}{E}$, i.e. we set $\Qh[E]{k}{d} = \Poly[d]{k}{E}$. 

In particular, our library offers some alternatives to properly define the degrees of freedom related to both the pressure and the velocity variables and to compute the polynomial projections of the virtual velocity basis functions. This list includes
\begin{itemize}
    \item The standard monomial approach detailed in \cite{DaVeiga2016}, where the sets $\GPoly{k}{\nabla}{E}$ and $\GPoly{k}{\perp}{E}$ are defined through the scalar monomials basis \eqref{eq:polynomials:setscaledmonomial}. The corresponding approach is implemented in the class \lstinline[language=C++, style=mystyle]{Polydim::VEM::MCC::VEM_MCC_2D_Velocity_LocalSpace} if $d=2$ and \lstinline[language=C++, style=mystyle]{Polydim::VEM::MCC::VEM_MCC_3D_Velocity_LocalSpace} if $d=3$.
    \item The partial andhonormal approaches introduced in \cite{Teora2024_mixed} by the authors. These approaches guarantee well-conditioned local and global matrices for the higher-values of the local polynomial degree $k$ with respect to the monomial approach. The related local spaces can be built exploiting \lstinline[language=C++, style=mystyle]{Polydim::VEM::MCC::VEM_MCC_2D_Partial_Velocity_LocalSpace} and \lstinline[language=C++, style=mystyle]{Polydim::VEM::PCC::VEM_MCC_2D_Ortho_Velocity_LocalSpace}, respectively, for $d=2$.
\end{itemize}
The aforementioned classes for the two-dimensional spaces use the Gauss quadrature points to define the boundary degrees of freedom. Orthonormal edge degrees of freedom are further available for the monomial and the orthonormal approaches. The construction of virtual projections through these degrees of freedom is detailed by the authors in \cite{Teora2023} and implemented in the classes \lstinline[language=C++, style=mystyle]{Polydim::VEM::MCC::VEM_MCC_2D_EdgeOrtho_Velocity_LocalSpace} and \lstinline[language=C++, style=mystyle]{Polydim::VEM::PCC::VEM_MCC_2D_Ortho_EdgeOrtho_Velocity_LocalSpace}.

\subsubsection{The Discrete Mixed Variational Formulation}

We define the global mixed virtual element spaces for both velocity and pressure variables as
\begin{gather*}
    \VMh{k}{}= \{\vv_h \in \VM \text{ s.t. } \vv_{h|E} \in \VMh[E]{k}{d}\ \forall E \in \Th\},\\
    \Qh{k}{} = \{q_h \in \Q \text{ s.t. } q_{h|E} \in \Qh[E]{k}{d}\ \forall E \in \Th\}.
\end{gather*}

As done for the primal formulation, we split the continuous bilinear form \eqref{eq:mixed:continuous_dbilin} according to the tessellation $\Th$ as
\begin{equation*}
\dbilin{\uu}{\vv} = \sum_{E \in \Th} \dbilin[E]{\uu}{\vv},\quad \dbilin[E]{\uu}{\vv} = \scal[E]{\K \uu}{ \vv}\quad \forall \uu,\vv \in \VM. 
\end{equation*}
Then, since we are not able to compute the quantity
\begin{equation}
\dbilin[E]{\uu_h}{\vv_h} = \scal[E]{\K \uu_h}{ \vv_h}\quad \forall \uu_h,\ \vv_h \in \VMh[E]{k}{d}, 
\label{eq:mixed:local_dbilin}
\end{equation}
we define the local virtual discrete bilinear form as the computable bilinear form
\begin{equation}
    \dbilinh[E]{\uu_h}{\vv_h} := \scal[0,E]{\K\, \vproj{0}{k} \uu_h}{ \vproj{0}{k} \vv_h} + \stab[E]{(I-\vproj{0}{k})\uu_h}{(I-\vproj{0}{k})\vv_h},
    \label{eq:mixed:ahE}
\end{equation}
where the stabilization term $\stab[E]{}{}$ is any symmetric and positive definite bilinear form that scales as \eqref{eq:mixed:local_dbilin}.
As in \cite{secondMixed,DaVeiga2016}, a trivial choice is represented by the dofi-dofi stabilization term:
\begin{equation*}
    \stab[E]{\uu_h}{\vv_h} = \vert E \vert \sum_{r=1}^{\Ndof[E]} \dof_r\left( \uu_h\right)\dof_r\left(\vv_h\right).
\end{equation*}

Finally, the mixed VEM approximation of \eqref{eq:mixed:variationalfomrulation} reads as: \textit{Find }$(\uu_{h}, p_h)  \in \VMh{k}{} \times \Qh{k}{} \textit{ such that}$:
\begin{equation}
    \begin{cases}
        \displaystyle\sum_{E\in\Th}\left[\dbilinh[E]{\uu_h}{\vv_h} -  \scal[E]{p_h}{ \nabla \cdot \vv_h}  \right] = 0 & \forall \vv_h \in \VMh{k}{}, \\
        \displaystyle\sum_{E\in\Th} \scal[E]{\nabla \cdot \uu_h}{ q_h}  = \displaystyle\sum_{E\in\Th} \scal[E]{f}{ q_h}  & \forall q_h \in  \Qh{k}{}.
    \end{cases}
    \label{eq:mixed:discreteformulation}
\end{equation}
The problem \eqref{eq:mixed:discreteformulation} has unique solution $(\uu_{h},p_h) \in \VMh{k}{} \times \Qh{k}{}$ and, for $h$ sufficiently small, the following a priori error estimates hold true
\begin{equation}
    \norm[\leb{2}{\Omega}]{p -p_h } = O( h^{k+1}),\quad 
    \norm[\vleb{2}{\Omega}{d}]{\uu - \uu_h} = O( h^{k+1}).
\end{equation}
Furthermore, the following superconvergence result holds true.
\begin{theorem}[Superconvergence result]
Let $p_h$ be the solution to \eqref{eq:mixed:discreteformulation} and let $ p_I \in \Qh{k}{}$ be the interpolant of $p$. Then, for $h$ sufficiently small, 
\begin{equation}
    \norm[\leb{2}{\Omega}]{ p_I - p_h} = O(h^{k+2}).
\end{equation}
\end{theorem}

\subsection{The Divergence-Free Formulation}\label{sec:DF}

Let us now consider the Stokes problem on $\Omega \subset \R^d$ with homogeneous Dirichlet boundary conditions:
\begin{equation}
    \begin{cases}
        \text{Find } (\uu,p) \text{ such that}\\
        - \nu \Delta \uu - \nabla p = \bm{f}& \text{in } \Omega,\\
        \div \uu = 0 & \text{in } \Omega,\\
        \uu = \bm{0} & \text{on } \Gamma = \partial \Omega,\\
    \end{cases}
    \label{eq:divfree:stokes_problem}
\end{equation}
where $\bm{f} \in \vleb{2}{\Omega}{d}$ represents the external source and $\nu \in \leb{\infty}{\Omega}$ is the uniformly positive viscosity.

Let us consider the spaces
\begin{equation}
    \VM := \vsob[0]{1}{\Omega}{d},\quad \Q := \leb[0]{2}{\Omega} = \left\{q \in \leb{2}{\Omega}: \int_{\Omega} q = 0\right\}.
\end{equation}
The variational formulation of the Stokes problem \eqref{eq:divfree:stokes_problem} reads as
\begin{equation}
    \begin{cases}
        \text{Find } (\uu,p) \in \VM \times \Q \text{ such that}\\
        \dbilin{\uu}{\vv} + \scal[\Omega]{  \div \vv}{p} = \scal[\Omega]{\bm{f}}{\vv} & \forall \vv \in \VM\\
        \scal[\Omega]{  \div \uu}{q} = 0 & \forall q\in \Q\\
    \end{cases}
    \label{eq:divfree:stokes_var_problem}
\end{equation}
where 
\begin{equation*}
    \dbilin{\uu}{\vv} = \int_{\Omega} \nu \nabla \uu : \nabla \vv \quad \forall \uu, \vv \in \VM \quad \text{ and } 
   \quad \nabla \uu = \begin{bmatrix}
       \frac{\partial u_1}{\partial x_1} & \cdots & \frac{\partial u_1}{\partial x_d} \\
       \vdots & \ddots & \vdots \\
      \frac{\partial u_d}{\partial x_1}  & \cdots & \frac{\partial u_d}{\partial x_d}
   \end{bmatrix}.
\end{equation*}

\subsubsection{The Local Space and the \lstinline[language=C++, style=mystyle]{Polydim::VEM::DF_PCC} classes}

Given $k \geq 2$, on each element $E \in \Th$, we define the following finite two-dimensional divergence-free virtual element space for the velocity variable $\uu$ \cite{Vacca2018}:
\begin{equation}
    \begin{aligned}
        \VMh[E]{k}{2} := \Big\{ \vv \in \vsob{1}{E}{2}:\ &(i)\ \ \ \vv_{|\partial E} \in \left[\Bk{k}{\partial E}\right]^2,\\
        &(ii)\ -\Delta \vv - \nabla s\in \GPoly{k}{\perp}{E} \text{ for some } s \in \leb[0]{2}{E},\\
        &(iii)\ \div \vv \in \Poly{k-1}{E} \\
        &(iv) \ \scal[E]{\vv - \proj{\nabla,E}{k} \vv}{\g^{\perp,k}} = 0\ \forall \g^{\perp,k} \in \GPoly{k}{\perp}{E} \setminus \GPoly{k-2}{\perp}{E} \Big\}.
\end{aligned}
\label{eq:divfree:discrete_local_velocity_space_stokes}
\end{equation}
If $d=3$, we define three-dimensional local primal divergence-free virtual element space for the velocity variable $\uu$ as
\begin{equation}
\begin{aligned}
\VMh[E]{k}{3} = \Big\{\vv \in \vsob{1}{E}{3}: &(i)\ \ \vv_{| \partial E} \!\in\! \left[\mathbb{U}_{k,1}(\partial E)\right]^3,\\
        &(ii)\ -\Delta \vv - \nabla s\in \GPoly{k}{\perp}{E} \text{ for some } s \in \leb[0]{2}{E},\\
        &(iii)\ \div \vv \in \Poly{k-1}{E} \\
        &(iv) \ \scal[E]{\vv - \proj{\nabla,E}{k} \vv}{\g^{\perp,k} } = 0\ \forall \g^{\perp,k}  \in \GPoly{k}{\perp}{E} \setminus \GPoly{k-2}{\perp}{E} \Big\}.
\end{aligned}
\label{eq:divfree:threelocalvirtualspace}
\end{equation}
where the space $\mathbb{U}_{k,1}(\partial E)$ is defined in \eqref{eq:primal:virtualfacespace}.
For the pressure variable $p$ we set
\begin{equation}
    \Qh[E]{k}{d} := \Poly[d]{k-1}{E}.
    \label{eq:divfree:discrete_local_pressure_space_stokes}
\end{equation}

Given a function $\vv \in \VMh[E]{k}{d}$, we consider the following degrees of freedom (see Figure \ref{fig:divfree:dofs} for a graphical illustration):
\begin{itemize}\label{dof:stokes_velocity}
    \item \textbf{Vertex DOFs}: the values of $\vv$ at the vertices of the polygon.
    \item \textbf{Edge DOFs}: the values of $\vv$ at $k-1$ Gauss-Lobatto internal points of every edge $e \in \Eh[E]$.
    \item \textbf{Face DOFs}: if $d=3$, the internal moments on each face $F \in \Fh[E]$, defined as
    \begin{equation}
        \frac{1}{\vert F \vert} \int_F \vv \cdot \nn_F\ \ppoly_{\alpha},\qquad \quad \forall I=1,\dots,2n^2_{k-2},\quad\frac{1}{\vert F \vert} \int_F \vv_{\bm{\tau}} \cdot \pp_I \quad \forall I=1,\dots,2n^2_{k-2},
    \end{equation}
    where $\vv_{\bm{\tau}} = \vv - \vv \cdot \nn_F$, $\{\pp_I\}_{I=1}^{2n_{k-2}^2}$ is a vector polynomial basis for $\VPoly[2]{k-2}{F}{2}$ built as in \eqref{eq:polynomials:vector_monomials} and $\{\ppoly_{\alpha}\}_{\alpha=1}^{n^2_{k-2}}$ is a scalar polynomial basis for $\Poly[2]{k-2}{F}$.
    \item \textbf{Internal $\bigoplus$ DOFs}: the moments of $\vv$
    \begin{equation}
        \int_E \vv \cdot \g^{\perp,k-2}_{\alpha} \qquad \forall \alpha=1,\dots,n^{\perp}_{k-2}.
    \end{equation}
    \item \textbf{Internal $\div$ DOFs}: the moments up to order $k-1$ and greater than zero of $\div \vv$ in $E$, i.e.
    \begin{equation}
        \int_E \div \vv \ \ppoly \qquad \forall \ppoly \in \mathcal{P}_{k-1}^d(E) \setminus \R.
    \end{equation}
\end{itemize}
It holds
\begin{gather}
    \Ndof[E] := \dim \VMh[E]{k}{d} = \begin{cases}
    2(\Nv[E] + \Ne[E] (k-1)) + n^{\perp}_{k-2} + n^2_{k-1} - 1& \text{if } d=2,\\
    3\left(\Nv[E] + \Ne[E] (k-1) + \Nf[E] \frac{(k-1)k}{2}\right) + n^{\perp}_{k-2} + n^3_{k-1} - 1  & \text{if } d=3,
    \end{cases}\\
    \dim \Qh[E]{k}{d} = \dim \Poly[d]{k-1}{E} = n_{k-1}^d.
\end{gather}

\begin{proposition}[\cite{Vacca2018, DaveigaDassi2020}]
    Let $\VMh[E]{k,\ell}{d}$ be the space defined in \eqref{eq:divfree:discrete_local_velocity_space_stokes} for $d=2$ and in \eqref{eq:divfree:threelocalvirtualspace} for $d=3$. The local DOFs, along with the enhancement property, allow us to compute exactly (up to machine precision) the face projections
    \begin{equation*}
        \proj{0,F}{k+1} : \left[\VPh[E]{k,1}{2}\right]^3 \to \VPoly[2]{k+1}{F}{3},\quad \proj{\nabla,F}{k} : \left[\VPh[E]{k,1}{2}\right]^3 \to \VPoly[2]{k}{F}{3}, 
    \end{equation*}
    and the element projections
    \begin{equation*}
        \proj{0,E}{k} : \VMh[E]{k}{d} \to \VPoly[d]{k-2}{E}{d},\quad \proj{\nabla,E}{k} : \VMh[E]{k}{d} \to \VPoly[d]{k}{E}{d}, \quad \vproj{0,E}{k-1} : \nabla \VMh[E]{k}{d} \to \VPoly[d]{k-1}{E}{d \times d}
    \end{equation*}
     as defined in \eqref{eq:proj:Pi0k} and \eqref{eq:proj:PiNabla}, i.e. we are able to compute the polynomials $\proj{0,F}{k+1} \vv_h$, $\proj{\nabla,F}{k} \vv_h$, $\proj{0,E}{k} \vv_h$, $\proj{\nabla,E}{k} \vv_h$ and $\proj{0}{k-1} \nabla \vv_h$, for any $\vv_h \in \VMh[E]{k}{d}$, using only the information given by the local degrees of freedom of $\vv_h$.
\end{proposition}

\subsubsection{The Discrete Divergence-Free Formulation}

We define the global virtual element spaces as
\begin{gather*}
    \VMh{k}{} := \{\vv \in \VM: \vv_{|E} \in \VMh[E]{k}{d} \ \forall E \in \Th \},\\
    \Qh{k}{} := \{q \in \Q: q_{|E} \in \Qh[E]{k}{d} \ \forall E \in \Th \}.
\end{gather*}
The divergence-free virtual discretization of Stokes problem reads as: $\textit{Find } (\uu_h,p_h) \in \VMh{k}{} \times \Qh{k}{} \textit{ such that}$
\begin{equation}
    \begin{cases}
        \sum_{E \in \Th}\left[\dbilinh[E]{\uu_h}{\vv_h} + \scal[E]{\div \vv_h}{p_h}\right] = \sum_{E \in \Th}\scal[E]{\bm{f}_h}{\proj{0,E}{k}\vv_h} & \forall \vv_h \in \VMh{k}{},\\
        \sum_{E \in \Th}\scal[E]{\div \uu_h}{q_h} = 0 & \forall q_h\in \Qh{k}{},\\
    \end{cases}
    \label{eq:divfree:stokes_discrete_problem}
\end{equation}
where we set
\begin{equation*}
    \dbilinh[E]{\uu_h}{\vv_h} = \scal[E]{\nu \proj{0,E}{k-1}\nabla \uu_h}{\proj{0,E}{k-1}\nabla \vv_h} + \stab[E]{(I - \proj{\nabla,E}{k})\uu_h}{(I - \proj{\nabla,E}{k})\vv_h}
\end{equation*}
and $\stab[E]{}{}$ is built following the construction detailed for the primal version of the method. If $(\uu_h,p_h)$ are the solution of \eqref{eq:divfree:stokes_discrete_problem}, it holds
\begin{equation*}
    \norm[\vleb{2}{E}{d\times d}]{\nabla \uu - \nabla \uu_h} = O(h^{k}),\quad \norm[\leb{2}{\Omega}]{p - p_h} =O(h^k).
\end{equation*}

We observe that this space offers the possibility to deal with more general problems, such as the Navier-Stokes equation and the Brinkman equation, being stable for both the Stokes and Darcy limit case \cite{DaVeigaLovadina2018, DaVigaMora2019, Vacca2018}.

\begin{remark}
    Let us consider that the velocity field $\uu$ should satisfy the continuous constraint $\div \uu = f$ for a certain scalar field $f$. We observe that
    \begin{equation}
        \div \VMh{k}{} = \Qh{k}{},
    \end{equation}
    and that, in particular, the divergence-free discrete solution $\uu_h$ satisfies \cite{Vacca2018}:
    \begin{equation}
        \div \uu_h = \proj{0,E}{k-1}f = \proj{0,E}{k-1} \div \uu.
    \end{equation}
    Thus, is $f=0$ as happens for the Stokes problem, $\uu_h$ is a ``point-wise'' divergence-free vector. On the other hand, employing the space detailed for the elastic problem in Section \ref{sec:elastic}, the solution 
    $\overline{\uu}_h$ of the corresponding discretization satisfies the divergence-free property only in a ``projected sense'', i.e.
    \begin{equation}
        \proj{0,E}{k-1} \div \overline{\uu}_h = \proj{0,E}{k-1}f = \proj{0,E}{k-1} \div \uu.
    \end{equation}
\end{remark}

\subsubsection{The Reduced Spaces}

In PolyDiM, there is the possibility to deal with the reduced version of the above method, which is introduced in \cite{DaVeigaLovadina2017}. In this reduced version, a large number of degrees of freedom
can be automatically eliminated from the system, obtaining one pressure DOF per element and a few internal DOFs for the velocity field. 
Precisely, the new spaces are given by
\begin{equation}
    \rVMh[E]{k}{d} = \{\vv \in \VMh[E]{k}{d}: \div \vv \in \Poly[d]{0}{E}\},\quad \rQh[E]{k}{d} = \Poly[d]{0}{E}.
\end{equation}
The subset given by the vertex, edge, (face) and the internal $\perp$ DOFs represents a unisolvent set of DOFs for the reduced velocity space $\rVMh[E]{k}{d}$, by leading to
\begin{gather*}
    \dim \rVMh[E]{k}{d} = \dim \VMh[E]{k}{d} - (n_{k-1}^d - 1) <  \dim \VMh[E]{k}{d},\\
    \dim \rQh[E]{k}{d} = 1 < \dim \Qh[E]{k}{d}.
\end{gather*}
\begin{proposition}[\cite{DaVeigaLovadina2017}]\label{prop:reduced_discrete_problem}
    Let $(\uu_h,p_h) \in \VMh{k}{} \times \Qh{k}{}$ be the solution of the discrete problem \eqref{eq:divfree:stokes_discrete_problem} and $(\widehat{\uu}_h,\widehat{p}_h) \in \rVMh{k}{} \times \rQh{k}{}$ be the solution of the related reduced discrete problem. Then
    \begin{equation}
        \widehat{\uu}_h = \uu_h \qquad \widehat{p}_{h|E} = \proj{0,E}{0}p_h \quad \forall E \in \Th.
    \end{equation}
\end{proposition}

Both the reduced and the full versions of the method are implemented in the library through classes belonging to the namespace \lstinline[language=C++, style=mystyle]{Polydim::VEM::DF_PCC}. We observe that, if $d=3$, these classes use local spaces defined in \lstinline[language=C++, style=mystyle]{Polydim::VEM::PCC} to evaluate boundary contributions.

\section{The Implementation Philosophy}\label{sec:code}

Our implementation approach adheres to the standard Galerkin methodology: we design classes that operate element-wise and then we glue together elemental contributions to assemble the global system associated with the PDE we aim to solve, as in a standard FEM solver. 

Nonetheless, as in a standard polytopal element method, these classes exploit element geometric properties, but operate independently of the values of these properties. Indeed, no element classification is performed based on element geometric properties, such as the number of vertices and edges (and faces), convexity or the presence of aligned edges and faces.

\subsection{The Polytope Geometry}\label{sec:polytope_geoemtry}
\begin{figure}[!h]
    \centering
    \begin{subfigure}{0.35\linewidth}
        \includegraphics[width=\linewidth]{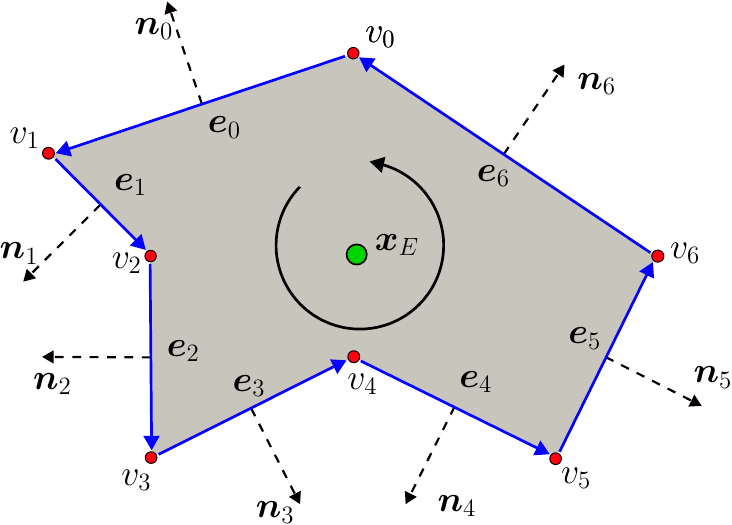}
        \caption{$d=2$.}
        \label{fig:geometry:2}
    \end{subfigure}\hspace{20pt}
    \begin{subfigure}{0.40\linewidth}
        \includegraphics[width=\linewidth]{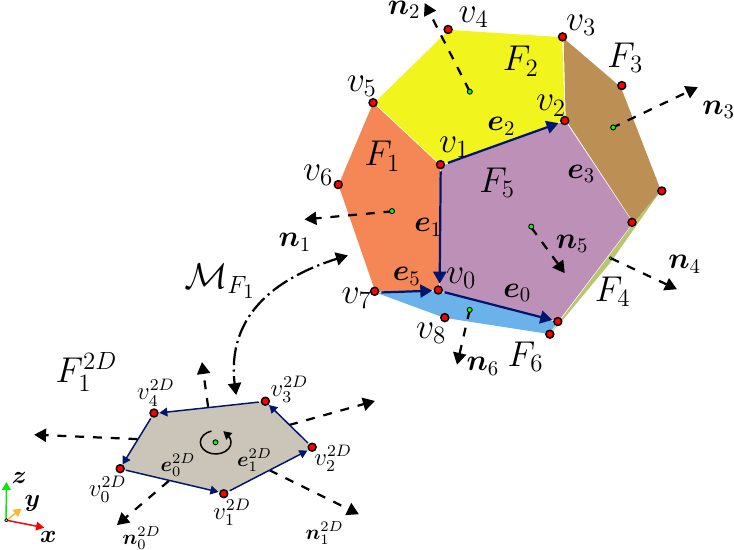}
        \caption{$d=3$.}
        \label{fig:geometry:3}
    \end{subfigure}
    \caption{Representation of Polygons and Polyhedrons.}
    \label{fig:geometry}
\end{figure}

Polygons are represented in PolyDiM as an ordered list of vertices coordinates. 
Specifically, for each polygon, we require that the list of its vertices is ordered counterclockwise, by implicitly defining a connectivity matrix for them. The order of edges is induced by the order of the vertices, as illustrated in Figure~\ref{fig:geometry:2}. For a polyhedron, instead, the required data include the list of vertices, along with the edges and faces connectivity matrices. The ordering of vertices within these connectivity matrices can be arbitrary. It is important to note that only the vertex ordering in the two-dimensional rotated face matters, since a polyhedron face also represents a polygon. Consequently, for each three-dimensional face $F$, a map $\mathcal{M}_{F}$ is needed to project the face onto the two-dimensional $(x_1,x_2,0)$-plane. After applying this transformation, the rotated face is represented as a two-dimensional polygon. Thus, we require that the map ensures the vertices of the resulting polygon are ordered counterclockwise. See Figure~\ref{fig:geometry:3} for an example.

Moreover, as in a standard Galerkin solver, the local discrete approximations may require other geometric properties.
For example, in Virtual Element we require the polytope diameter, centroid, measure and unit outward normals on each edge or face. All these geometric properties can either be supplied by the user or automatically computed using GeDiM.

\subsection{Integration}\label{sec:integration}
\begin{figure}[!h]
    \centering
    \begin{subfigure}{0.30\linewidth}
        \includegraphics[width=0.9\linewidth]{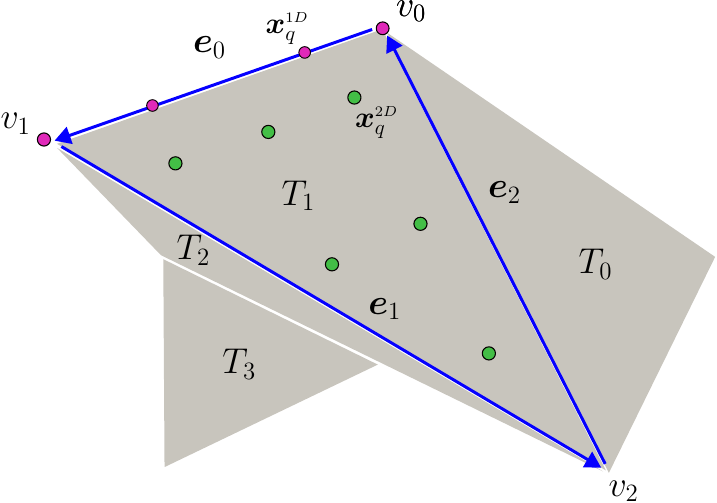}
        \caption{$d=2$.}
        \label{fig:quadrature:2}
    \end{subfigure}
    \begin{subfigure}{0.40\linewidth}
        \includegraphics[width=0.8\linewidth]{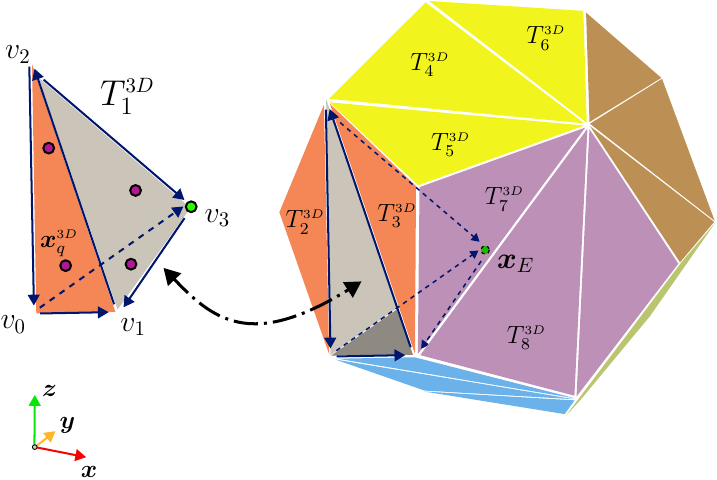}
        \caption{$d=3$.}
        \label{fig:quadrature:3}
    \end{subfigure}
    \caption{Sub-triangulation for a polygon (Left) and a polyhedron (right). green dots represent quadrature points related to a rule of order $4$ over a triangle, whereas red dots represent points corresponding to a quadrature formula of order $2$. }
    \label{fig:quadrature}
\end{figure}

Since PolyDiM is designed to deal with variational problems, we need to compute integrals. Moreover, to deal with polytopal methods, the integrals must be computed over generic polytopes. 
Integration of a generic polytope can be handled either by sub-triangulating the polytope and then resorting to a standard reduced Gaussian quadrature formula over triangles or tetrahedra, or by exploiting quadrature rules built ad-hoc for general polytope. For the sake of generality, we adopt the standard technique and proceed by sub-triangulating the polytope. To limit the computational cost due to the high density of the induced quadrature points, we make available in GeDiM many efficient algorithms to triangulate polytopes. Clearly, as with the geometric properties, both the sub-triangles and sub-tetrahedra can also be provided by the user.

Figures~\ref{fig:quadrature} illustrate an example of such polygon and polyhedron subdivisions. The sub-triangulation shown in Figure~\ref{fig:quadrature:2} is generated using the ear clipping algorithm provided by the GeDiM library. Note that the 2D triangles must preserve a counterclockwise ordering of their vertices.
In contrast, the sub-tetrahedra depicted in Figure~\ref{fig:quadrature:3} are constructed by connecting the triangulations of the polyhedron faces to its centroid. 
Alternatives to the methods described above are also available in the GeDiM library for both polygons and polyhedra.

Finally, edge and face integrals must also be computed, both to enforce Neumann boundary conditions and to account for boundary contributions in the VEM projection computations. Face integrals are evaluated using aforementioned polygonal integration techniques, whereas, for edge integrals, both Gauss and Gauss-Lobatto quadrature rules are available.

\subsection{Assembling}
\label{sec:assemble}
\begin{figure}[!h]
    \centering
    \begin{subfigure}{0.40\linewidth}
        \includegraphics[width=\linewidth]{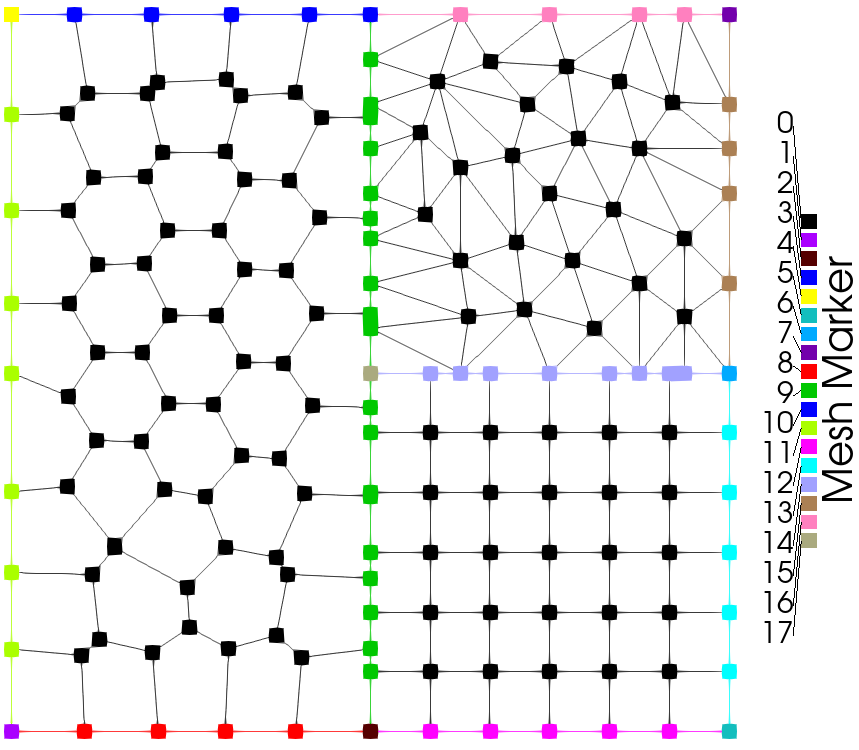}
        \caption{}
        \label{fig:marker:mesh}
    \end{subfigure}
    \begin{subfigure}{0.40\linewidth}
        \includegraphics[width=\linewidth]{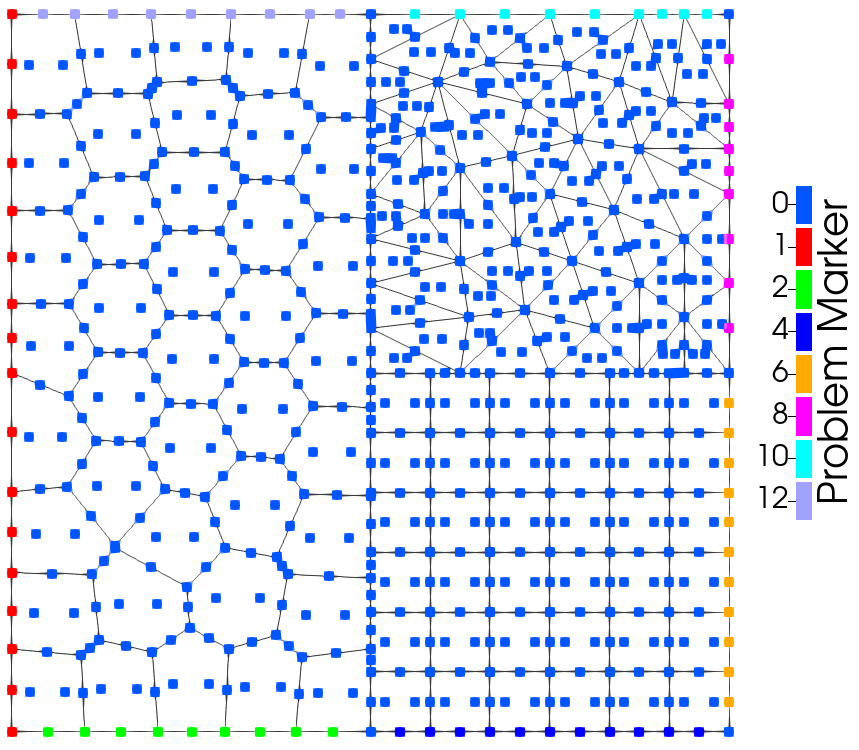}
        \caption{}
        \label{fig:marker:dofs}
    \end{subfigure}
    \caption{Left: Mesh cells coloured by mesh markers. Right: DOF cells coloured by problem markers.}
    \label{fig:marker}
\end{figure}
\begin{sflisting}[language=C++,style=mystyle,caption={Vendermonde matrices containing the projection of basis function derivatives evaluated at quadrature points.},label={lst:basis_der}]
    std::vector<Eigen::MatrixXd> vandermonde_derivatives = ComputeBasisFunctionsDerivativeValues(
        local_space_data,
        Polydim::VEM::PCC::ProjectionTypes::Pi0km1Der
    );
\end{sflisting}
\begin{lstlisting}[language=C++,style=mystyle,caption={Dofi-Dofi VEM stabilization matrix.},label={lst:stab}]
    Eigen::MatrixXd S = ComputeDofiDofiStabilizationMatrix(
        local_space_data,
        Polydim::VEM::PCC::ProjectionTypes::PiNabla
    );
\end{lstlisting}

The assembly workflow for all the implemented discretization methods in PolyDiM is based on the \emph{global-to-local} scheme typical of the Finite Element Method \cite{brennerscott}. Thus, the assembly process of a generic problem is divided into two main steps:
\begin{enumerate}
    \item Compute local matrices that represent each bilinear or linear operators dictated by the PDE.
    \item Assemble these local matrices into the global one.
\end{enumerate}

Specifically, all the implemented local spaces, including the ones detailed in Sections \ref{sec:PCC}, \ref{sec:MCC} and \ref{sec:DF}, are designed to evaluate the local (projected for the VEM) basis functions $\{\varphi_i\}_{i=1}^{\Ndof[E]}$ and their derivatives at points belonging to the element $E$. These values are then used to compute local matrices. 
For instance, the VEM elemental stiffness matrix $\mathbf{A} \in \R^{\Ndof[E] \times \Ndof[E]}$ related to the primal formulation presented in Section \ref{sec:PCC:discreteformulation} can be computed as:
\begin{align*}
    \mathbf{A}^E_{ij} &= \int_{E} \D \proj{0,E}{k-1} \nabla \varphi_j \cdot \proj{0,E}{k-1} \varphi_i + \stab[E]{(I-\proj{\nabla,E}{k}) \varphi_i}{(I-\proj{\nabla,E}{k}) \varphi_j} \\
    &\approx \sum_{q=1}^{N^Q_E} \omega_q \D(\xx_q) \proj{0,E}{k-1} \nabla \varphi_j(\xx_q) \cdot \proj{0,E}{k-1} \varphi_i (\xx_q) + \mathbf{S}^E_{ij} \quad \forall i,j=1,\dots,\Ndof[E],
\end{align*}
where $\{(\omega_q,\xx_q)\}_{q=1}^{N^Q_E}$ is a quadrature formula over $E$ with $N^Q_E$ points.
In this case, the class method of Code~\ref{lst:basis_der} returns the projection values of the gradients of virtual basis functions at the quadrature points of the element $E$ stored in Vandermonde matrices, whereas the Code~\ref{lst:stab} allows the user to retrieve the dofi-dofi stabilization matrix $\mathbf{S}^E$, defined in Equation \eqref{eq:primal:dofidofi}.

Concerning the second step, the global assembly phase inserts the entries of the local matrices into the global matrix associated with the discrete problem.
This process is performed through a map that associates each local degree of freedom of the local spaces to the corresponding global degree of freedom.
An example of the implementation of this \emph{local-to-global} map is contained inside the library class \lstinline[language=C++, style=mystyle]{Polydim::PDETools::DOFs::DOFsManager} and it is used in all the examples contained in the library.

We observe that, based on the formulation employed, one or more instances of \lstinline[language=C++, style=mystyle]{DOFsManager} can be used. For instance, to solve the mixed differential problem \eqref{eq:mixed:variationalfomrulation} we employ two instances of \lstinline[language=C++, style=mystyle]{DOFsManager}. 
The first counts the degrees of freedom associated with the velocity variable, while the second accounts for pressure DOFs. 
Similarly, to deal with the elasticity problem \eqref{eq:primal:vem_elasticity} we employ again two objects of type \lstinline[language=C++, style=mystyle]{DOFsManager}. In this case, each is associated with a different component of the discrete (vector)-solution.

The local-to-global map is constructed based on the local discrete space information contained within the \lstinline[language=C++, style=mystyle]{ReferenceElement} classes. Since the VEM can handle meshes composed of arbitrarily shaped polygons and polyhedra, the concept of reference element in this context differs from the geometric interpretation typically used in FEM. Indeed, this {ReferenceElement} serves as a container for local information that is common across all elements in the tessellation: the number of degrees of freedom associated with vertices (called DOF cells $0$D), edges (DOF cells $1$D), faces (DOF cells $2$D) if $d=3$, and internal to the element $E$ (DOF cells $d$D), the order of quadrature rules employed to compute projections and the information about objects related to the polynomial basis used to define both DOFs and projections.

Moreover, this map helps to store information about the test problem and, in particular, about the type of conditions imposed on each considered elemental DOF cell. Indeed, to impose different conditions dictated by the problem (e.g. boundary and interface conditions) on different parts of the domain, each geometric cell of the mesh $\Th$ is marked by a number. Figure~\ref{fig:marker:mesh} shows an example of these \emph{mesh markers} in the case of a domain divided into three sub-domains according to problem coefficients values: marker $9,\ 13$, and $17$ denote geometric cells belonging to the sub-domains interfaces, whereas all the other markers denotes different parts of the boundary domains. The only exception in this example is represented by the marker $0$ that labels mesh cells that does not correspond to the particular condition. 
Particularly, all the $d$D geometric cells are always marked by $0$. 

Given a set of mesh markers, the test problem associates with each specific mesh marker a \emph{problem marker}, which provides information about the boundary condition that must be applied to the DOF cell related to this geometric cell.
Given this information, the {DOFsManager} tells us if a local DOF cell corresponds to
\begin{itemize}
    \item a \emph{strong} cell, i.e. an essential condition which is imposed strongly in the discrete formulation \cite{brennerscott}. This is the case, for instance, of cells belonging to the Dirichlet boundary for the primal version of the method or to the Neumann boundary for the mixed formulation.
    \item a \emph{weak} cell, i.e. a natural condition which is imposed weakly in the discrete formulation \cite{brennerscott}. The weak cells correspond to DOFs (unknown) of the discrete problem. This is the case, for instance, of cells belonging to the Neumann boundary for the primal version of the method or to the Dirichlet boundary for the velocity variable in the mixed formulation.
    \item a cell with where no boundary condition is imposed that correspond to DOFs for the discrete problem. This can be the case of internal cells or $(d-2)$D cells in between two weak cells.
\end{itemize}

An example of this process is shown in Figure~\ref{fig:marker:dofs} for the DOF cells related to the virtual divergence-free velocity space \eqref{eq:divfree:discrete_local_velocity_space_stokes} for $k=2$. This figure illustrates how the {DOFsManager} associates with each DOF cell a problem marker. In this example, the problem marker $0$ denotes DOFs for which no boundary conditions are imposed. In contrast, odd problem markers greater than 
$0$ correspond to strong DOF cells, whereas even problem markers greater than 
$0$ correspond to weak DOF cells.

\subsection{Solver}
\label{sec:solver}
After the assembly phase, to find the discrete solution, i.e. the vector of DOFs, we must solve the global system of equations, that can be a linear or non-linear systems of equations. For this purpose, we offer through GeDiM an interface with the most common linear algebra C++ libraries, including \texttt{Eigen} \cite{eigenweb}, and \texttt{PETSc} \cite{petsc-web-page}.

%For example, users can refer to the Cholesky decomposition interface for \texttt{Eigen}, implemented in the class \lstinline[language=C++, style=mystyle]{Gedim::Eigen_CholeskySolver}, or to the interface for PETSc's KSP linear solvers, available in the class \lstinline[language=C++, style=mystyle]{Gedim::PETSc_KSPSolver}.
%It is worth noting that all example applications provided in the library are implemented using \texttt{Eigen} solvers, to avoid requiring users to install \texttt{PETSc} on their machines.

All solver interfaces are designed to operate with the generic classes \lstinline[language=C++, style=mystyle]{Gedim::ISparseArray} and \lstinline[language=C++, style=mystyle]{Gedim::IArray}, which represent the sparse matrices and vectors of the linear system, respectively.
Concrete implementations of these abstract classes are already available for both \texttt{Eigen} and the serial version of \texttt{PETSc} within the GeDiM library.
This level of abstraction in the linear algebra layer, if used by the user, enables to switch between different solvers with ease and without altering the core components of the code, such as the assembly process.

\section{Usage and Examples}\label{sec:examples}

The numerical solution of partial differential equations on geometrically complex domains represents a fundamental task in scientific computing, owing to its central role in accurately modeling a broad spectrum of real-world phenomena, but it can be extremely challenging. Polytopal methods like the Virtual Element Method alleviate these meshing issues by allowing, for instance, \cite{Attene2021}
\begin{itemize}
    \item to achieve great flexibility in the treatment of complex geometries; 
    \item automatic handling of hanging nodes;
    \item simplify refinement and agglomeration strategies.
\end{itemize}

In this section, we present a selection of numerical experiments that demonstrate the flexibility and robustness of our library in handling complex geometries. Most of the examples are drawn from the authors previous research, where PolyDiM has been effectively employed to solve challenging problems using advanced polytopal methods. We briefly summarize these applications to highlight the capabilities of the library, while referring the reader to the original papers for comprehensive details. Additionally, we include examples from external studies where PolyDiM could be successfully applied, further illustrating its versatility and broad applicability across various scientific computing contexts.
\begin{figure}[!ht]
    \centering
    \begin{subfigure}{0.33\linewidth}
        \includegraphics[width=\linewidth]{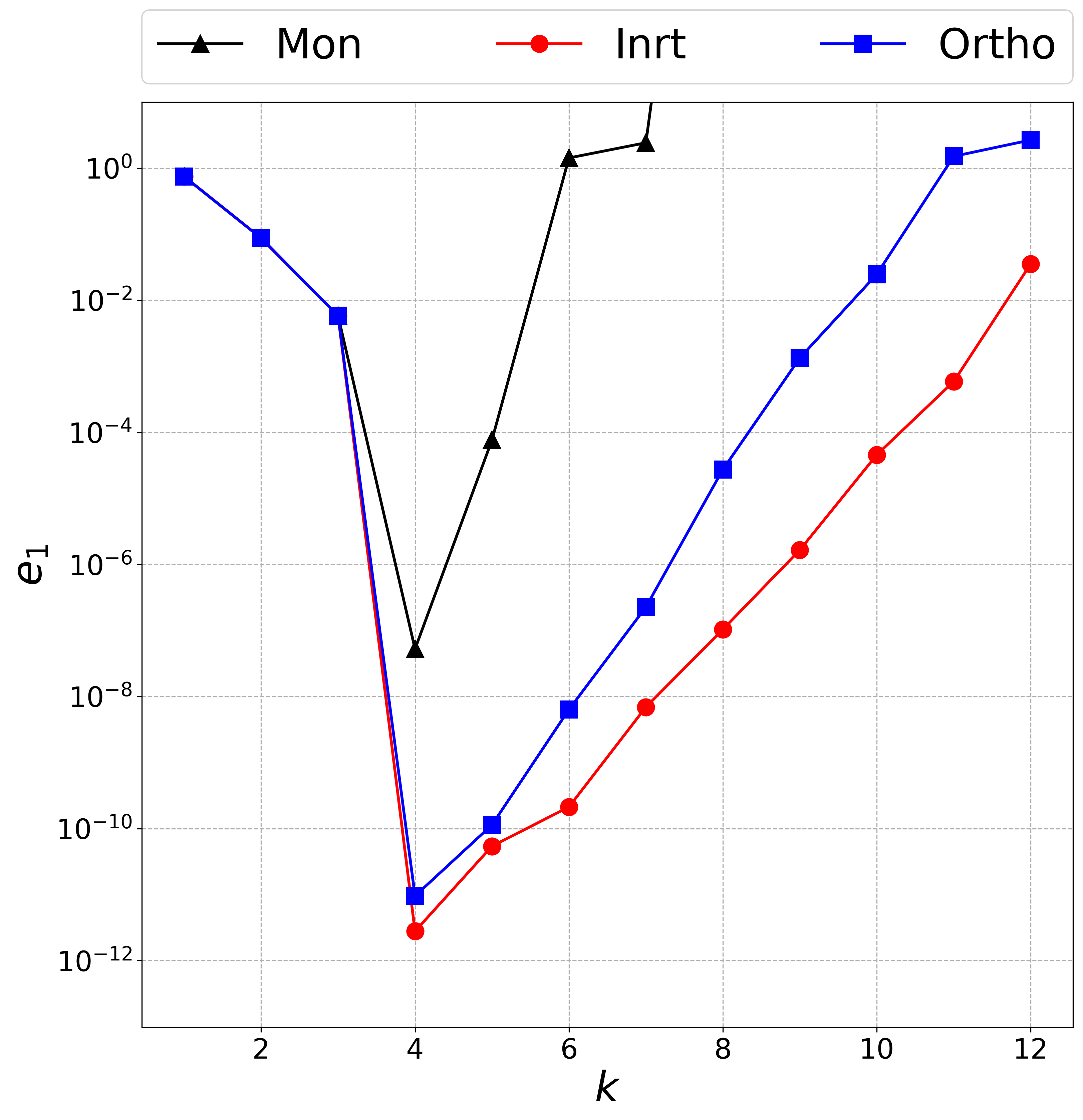}
        \caption{Test 1 in \cite{Teora2024}.}
    \end{subfigure}
    \begin{subfigure}{0.33\linewidth}
        \includegraphics[width=\linewidth]{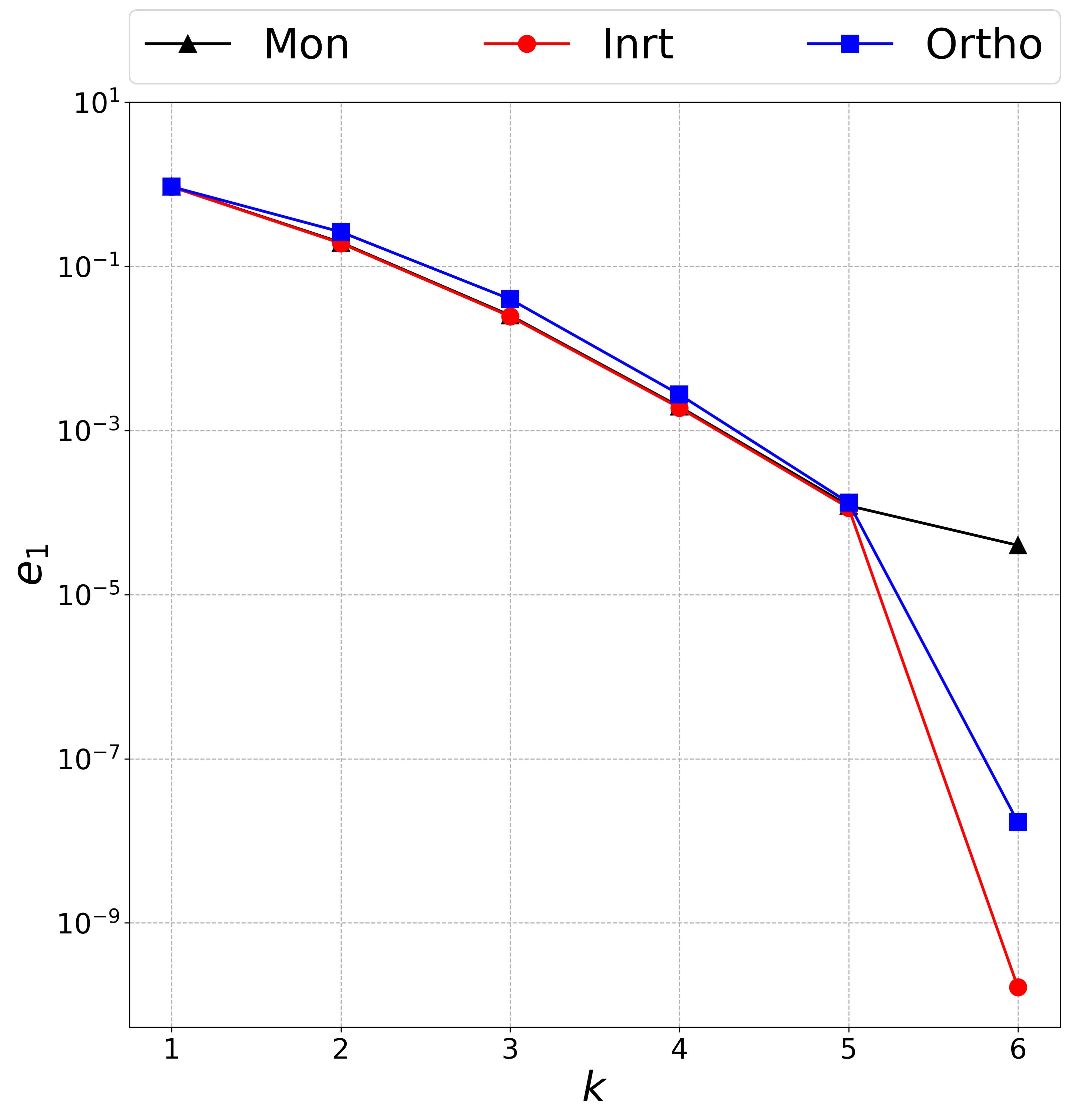}
        \caption{Test 3 in \cite{Teora2024}.}
    \end{subfigure}
    \begin{subfigure}{0.33\linewidth}
        \includegraphics[width=\linewidth]{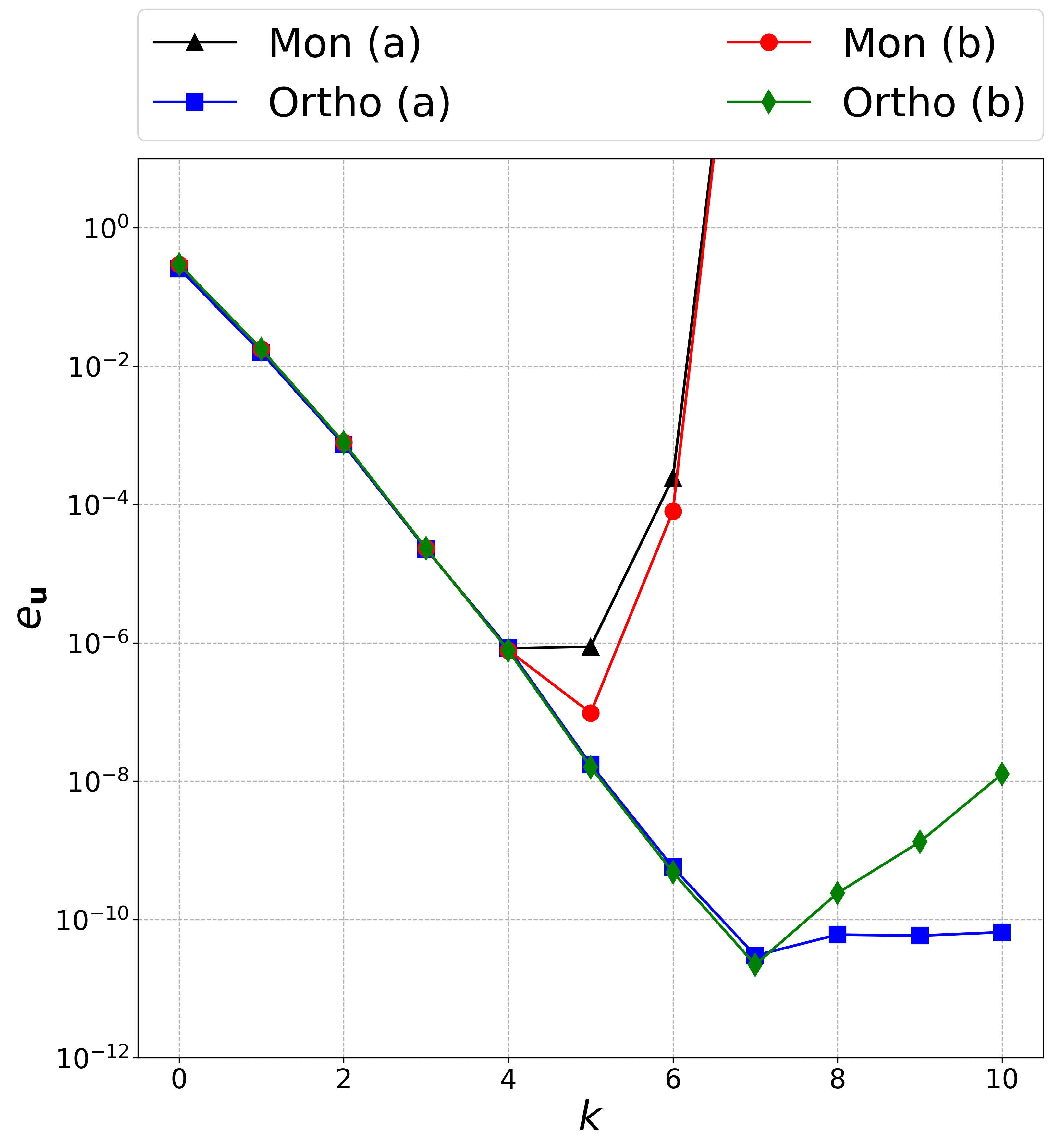}
        \caption{Test 1 in \cite{Teora2023}.}
    \end{subfigure}
    \caption{Left: $H^1$-error $e_1$ for the primal version of the method as the method order $k$ varies in two-dimensions for the different approaches proposed in Section \ref{sec:PCC}. Center: $H^1$-error for the primal version of the method as $k$ varies in three-dimensions for the different approaches proposed in Section \ref{sec:PCC}. Right: $L^2$-error for the velocity variable in the mixed version of the method as $k$ varies in two-dimensions for different approaches proposed in Section \ref{sec:MCC}. }
    \label{fig:test1:error}
\end{figure}

\subsection{Elliptic problems and Polynomial bases}

\begin{figure}[!ht]
    \centering
    \begin{subfigure}{0.33\linewidth}
        \includegraphics[width=\linewidth]{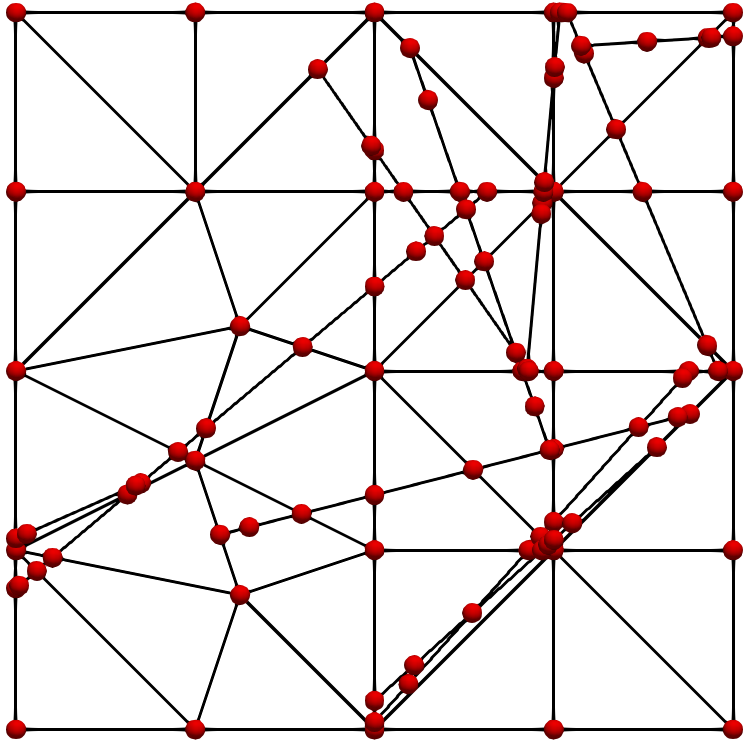}
        \caption{Test 1 in \cite{Teora2024}.}
    \end{subfigure}
    \begin{subfigure}{0.33\linewidth}
        \includegraphics[width=\linewidth]{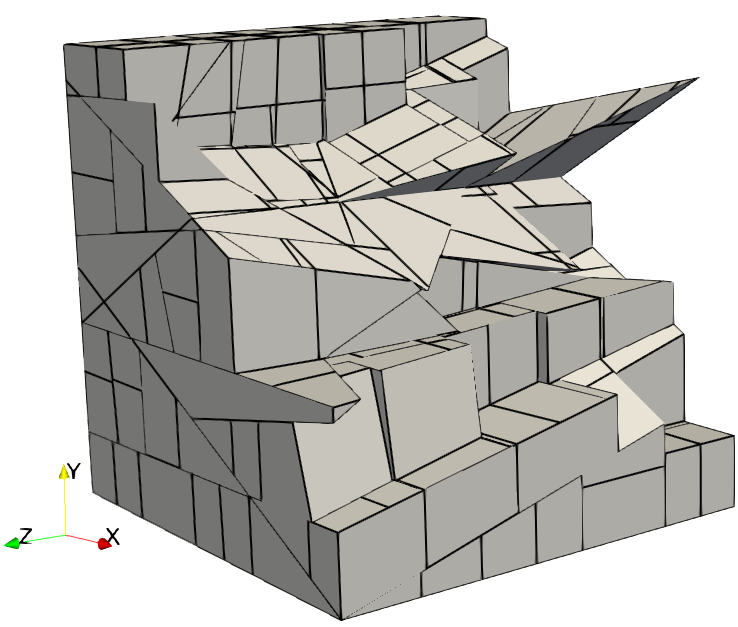}
        \caption{Test 3 in \cite{Teora2024}.}
    \end{subfigure}
    \begin{subfigure}{0.33\linewidth}
        \includegraphics[width=\linewidth]{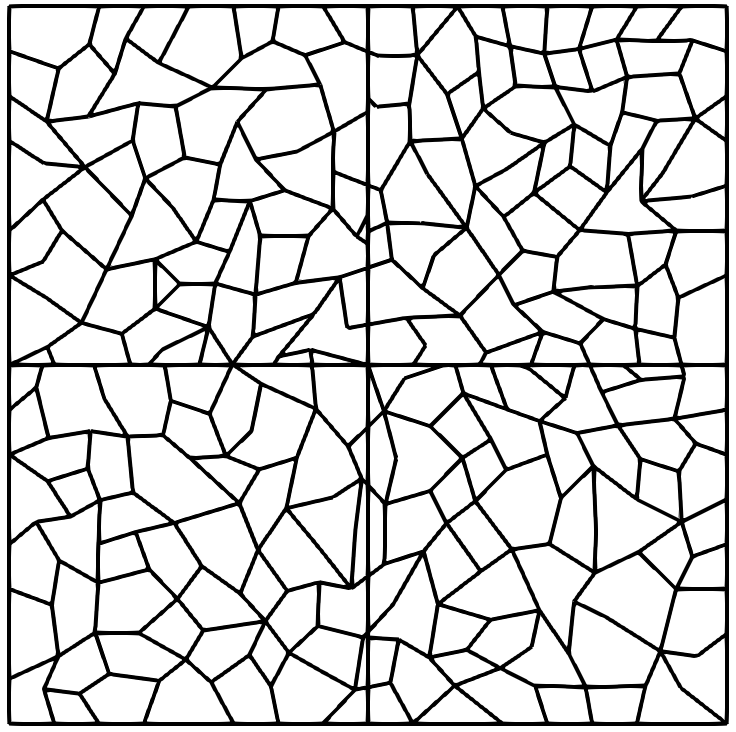}
        \caption{Test 1 in \cite{Teora2023}.}
    \end{subfigure}
    \caption{Meshes characterized by concave elements and a copious number of hanging nodes taken from Discrete Matrix and Fractures applications.}
    \label{fig:test1:mesh}
\end{figure}

\begin{figure}[!ht]
    \centering
    \begin{subfigure}{0.35\linewidth}
        \includegraphics[width=0.9\linewidth]{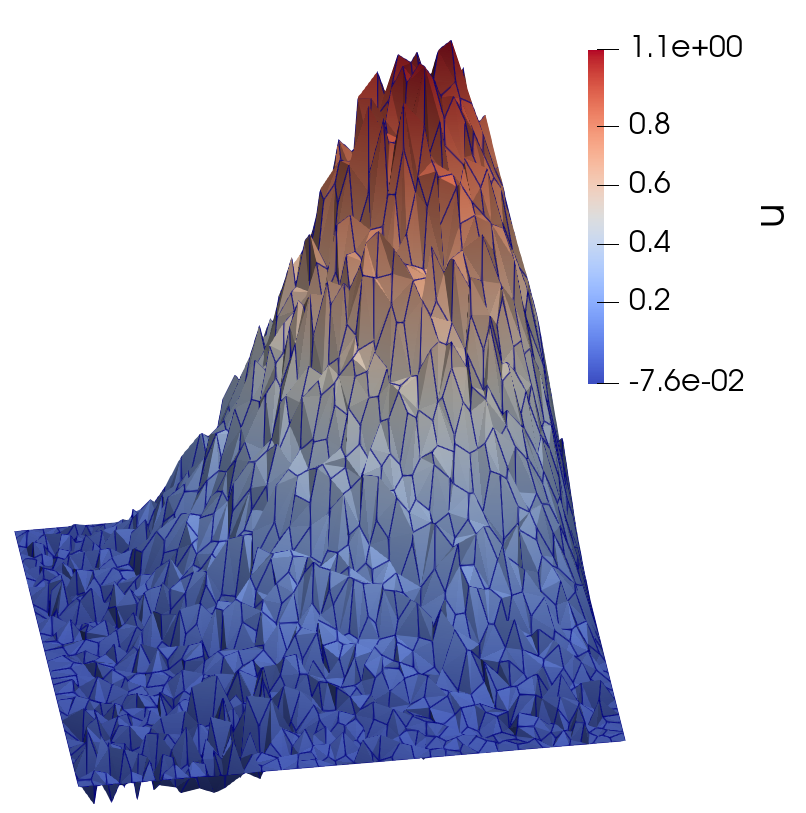}
        \caption{}
    \end{subfigure}
    \begin{subfigure}{0.35\linewidth}
        \includegraphics[width=0.9\linewidth]{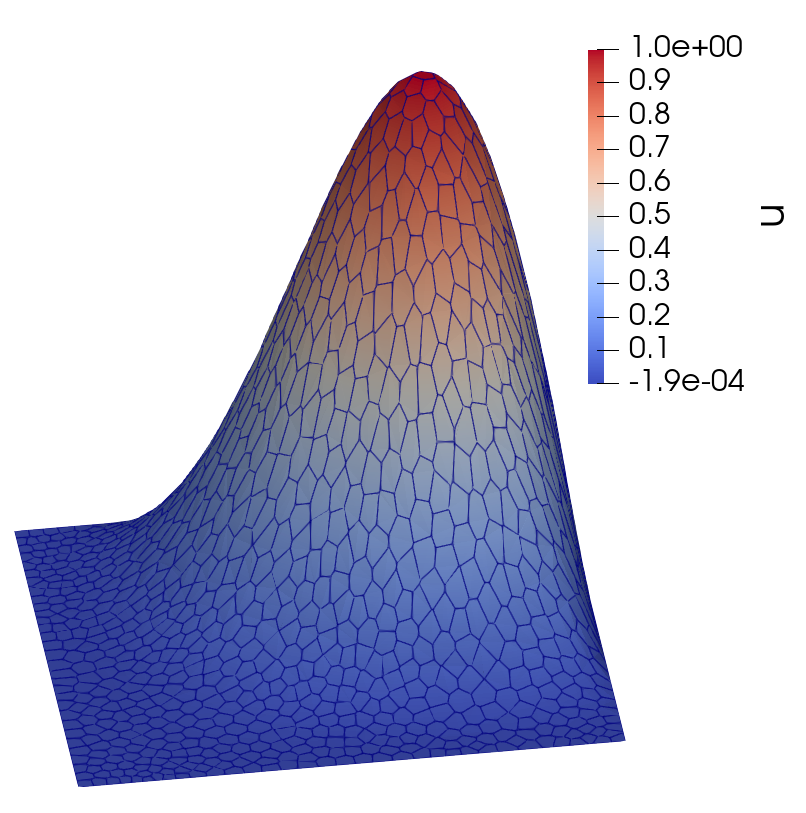}
        \caption{}
    \end{subfigure}
    \caption{A SUPG-stabilized VEM solution related to advection-dominated problem. Test described in \cite{Benedetto2016SUPG}.}
    \label{fig:test2:solution}
\end{figure}

PolyDiM enables the solution of general second-order elliptic problems on both two- and three-dimensional meshes, for a generic order $k$. Moreover, as discussed in the previous sections, PolyDiM provides methods that leverage novel polynomial bases to improve the conditioning of the discrete systems in the presence of badly-shaped polytope and for high orders of the method. Figure \ref{fig:test1:error} illustrates the behaviour of VEM errors as the order $k$ increases, when solving the general second-order elliptic problem on meshes shown in Figure \ref{fig:test1:mesh} for the different approaches detailed in Section \ref{sec:PCC}. These meshes are characterized by concave elements and a copious number of hanging nodes, representing typical examples of meshes encountered in Discrete Matrix and Fractures applications. These examples are taken from \cite{Teora2023, Teora2024, Teora2024_mixed}, where the authors have effectively used the PolyDiM library. Moreover, we recall that authors also contribute to introduce a consistent Streamline Upwind Petrov-Galerkin (SUPG) formulation for the virtual element methods. Figure shows the results related to Test 1 performed in \cite{Benedetto2016SUPG}, where an advection-dominated problem is solved. All these examples are implemented and available to the reader on GitHub.

\subsection{Elasticity, Brinkman and Navier-Stokes Problems}

\begin{figure}[!ht]
    \centering
    \includegraphics[width=0.3\linewidth]{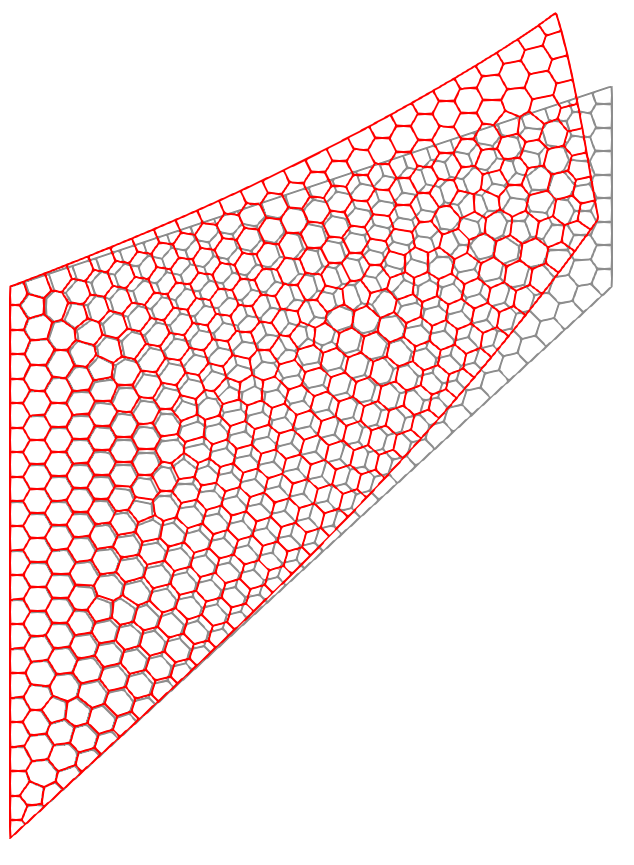}
    \caption{Test 2: Cook’s membrane. Undeformed (gray) and deformed (red) configurations. $k=1$.}
    \label{fig:test3:displacement}
\end{figure}

\begin{figure}[!ht]
    \centering
    \begin{subfigure}{0.33\linewidth}
        \includegraphics[width=\linewidth]{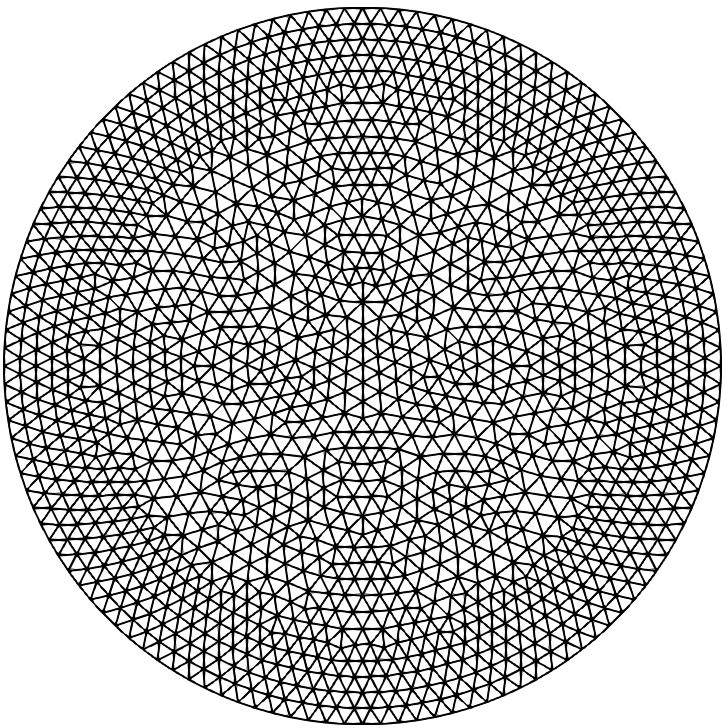}
        \caption{}
    \end{subfigure}
    \begin{subfigure}{0.33\linewidth}
        \includegraphics[width=\linewidth]{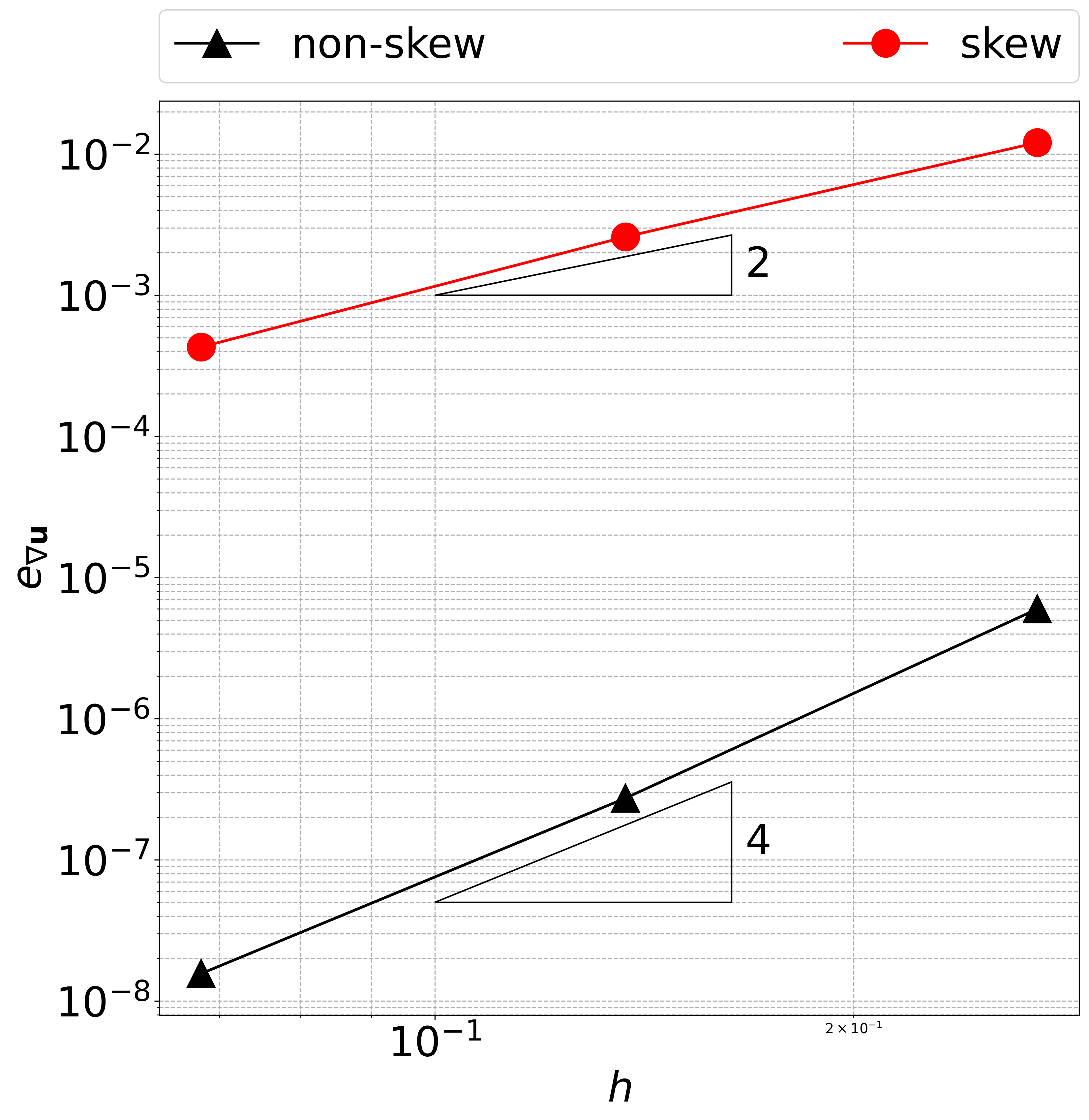}
        \caption{}
    \end{subfigure}
    \begin{subfigure}{0.33\linewidth}
        \includegraphics[width=\linewidth]{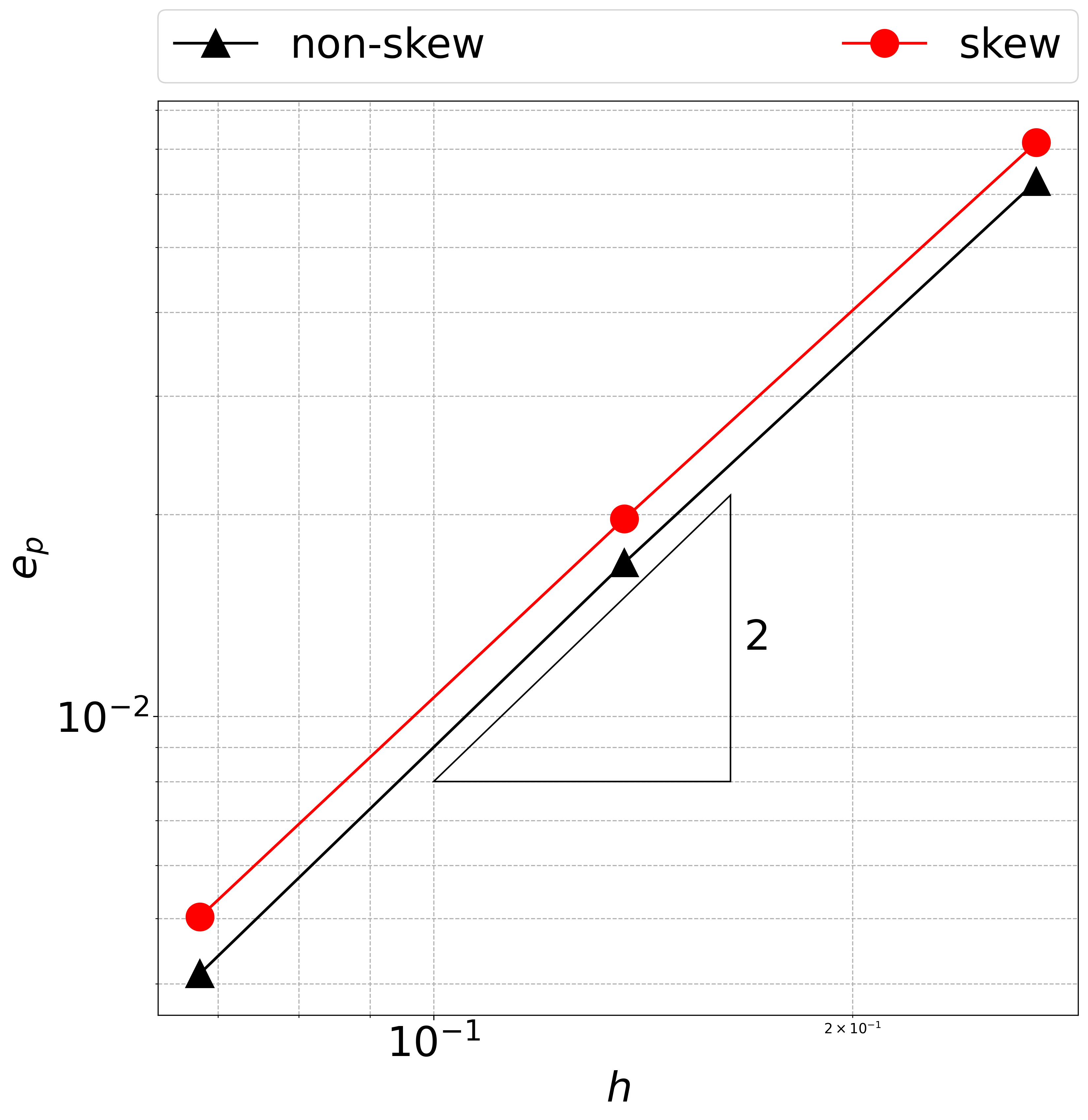}
        \caption{}
    \end{subfigure}
    \caption{Left: Domain and the finest mesh used to solve the Navier-Stokes problem. Center: Behaviour of errors related to the velocity variable as the mesh parameter $h$ decreases for both the skew and the non-skew trilinear form. Right: Behaviour of errors related to pressure variable as the mesh as the mesh parameter $h$ decreases for both the skew and the non-skew trilinear form. $k=2$.}
    \label{fig:test8:error}
\end{figure}

\begin{figure}[!ht]
    \centering
    \begin{subfigure}{0.4\linewidth}
        \includegraphics[width=1\linewidth]{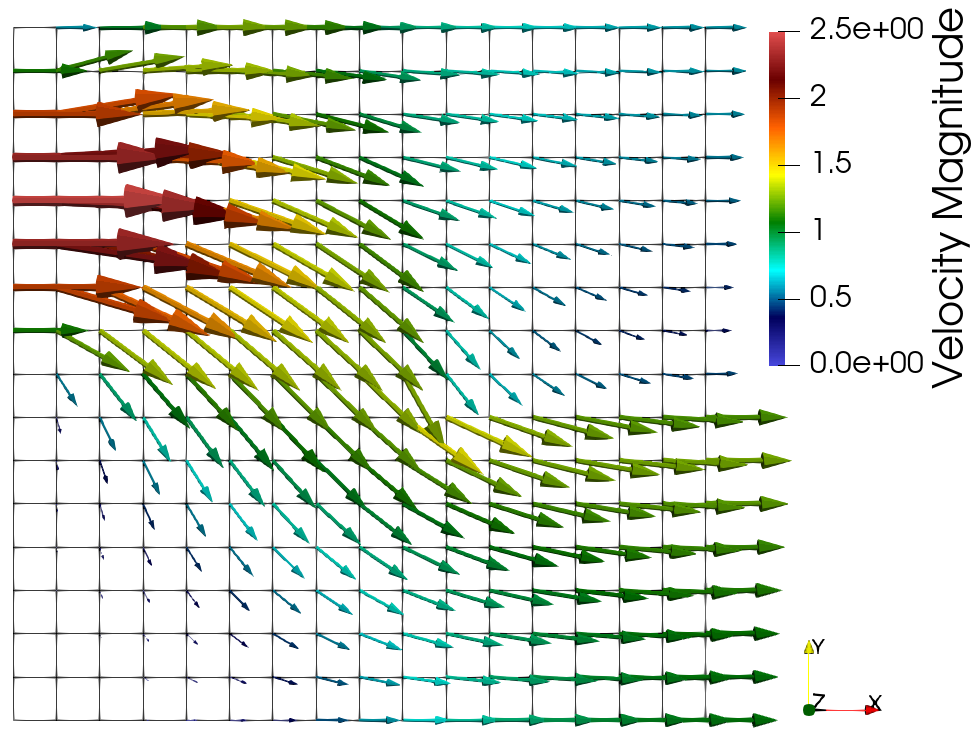}
        \caption{}
    \end{subfigure}\hspace{2pt}
    \begin{subfigure}{0.5\linewidth}
        \includegraphics[width=\linewidth]{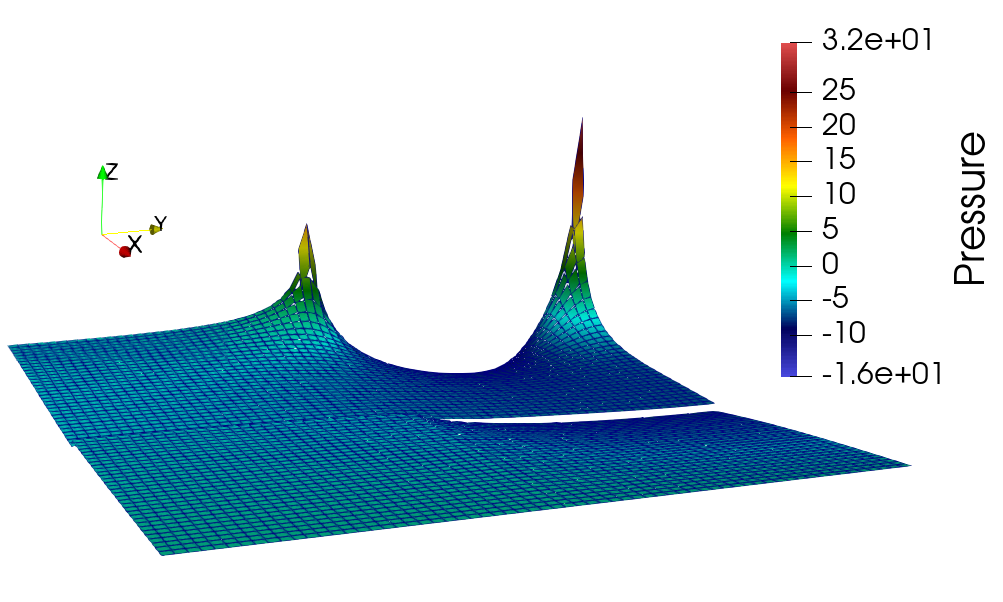}
        \caption{}
    \end{subfigure}
    \caption{Left: Velocity field computed by solving the Brinkman equation related to a coarse square mesh obtained with divergence-free formulation of virtual element for $k=2$. Right: The pressure field obtained related to a fine square mesh obtained with divergence-free formulation of virtual element for $k=2$.}
    \label{fig:test7:solution}
\end{figure}

PolyDiM supports the solution of various types of PDEs, including Brinkman, elasticity, and Navier–Stokes problems. In particular, divergence-free virtual elements can be effectively used to solve both the Navier–Stokes and Brinkman equations, offering stability in both the Darcy and Stokes limits (see Section \ref{sec:DF}), as well as in the incompressible limit for elasticity problems. Elasticity problems can also be addressed using standard primal formulations under compressibility assumptions, as detailed in Section \ref{sec:elastic}. The library provides tools for solving all these problems in both two- and three-dimensional settings.

To showcase the capabilities of PolyDiM, we provide a representative numerical experiment for each problem type. All examples are fully implemented and publicly available on GitHub.

For the elasticity case, inspired by Test 4 in \cite{Artioli2017}, we solve the $2D$ Cook’s membrane problem \eqref{eq:primal:vem_elasticity} with 
$k=1$ on a Voronoi mesh. We highlight that standard Finite Element Methods on triangular meshes are also available in PolyDiM. Since the exact solution is not provided for this example, a reference solution can be computed through this standard method to test VEM accuracy. 

Concerning fluid dynamic problems, Figure \ref{fig:test8:error} reports the VEM error corresponding to the resolution of the Navier–Stokes problem described in Test 5.2 of \cite{Artioli2017}, whereas Figure \ref{fig:test7:solution} illustrates the velocity and pressure fields obtained solving a Brinkman problem inspired by Test 5.6 in \cite{Artioli2017}. Notably, in this latter test, the DOF cells related to the boundary degrees for the first component of the velocity are associated with the problem markers as shown in Figure \ref{fig:marker:dofs}. In addition, the domain subdivision highlighted in the figure allows us to set different coefficient values across subregions. This setup enables us the simultaneous resolution of a Stokes problem in the left part of the domain and a pure Darcy problem in the right part, effectively illustrating the versatility of PolyDiM in handling mixed regimes within a single simulation.

\subsection{Mesh Quality Improvement}
\begin{figure}[!h]
    \centering
    \begin{subfigure}{0.35\linewidth}
        \includegraphics[width=0.9\linewidth]{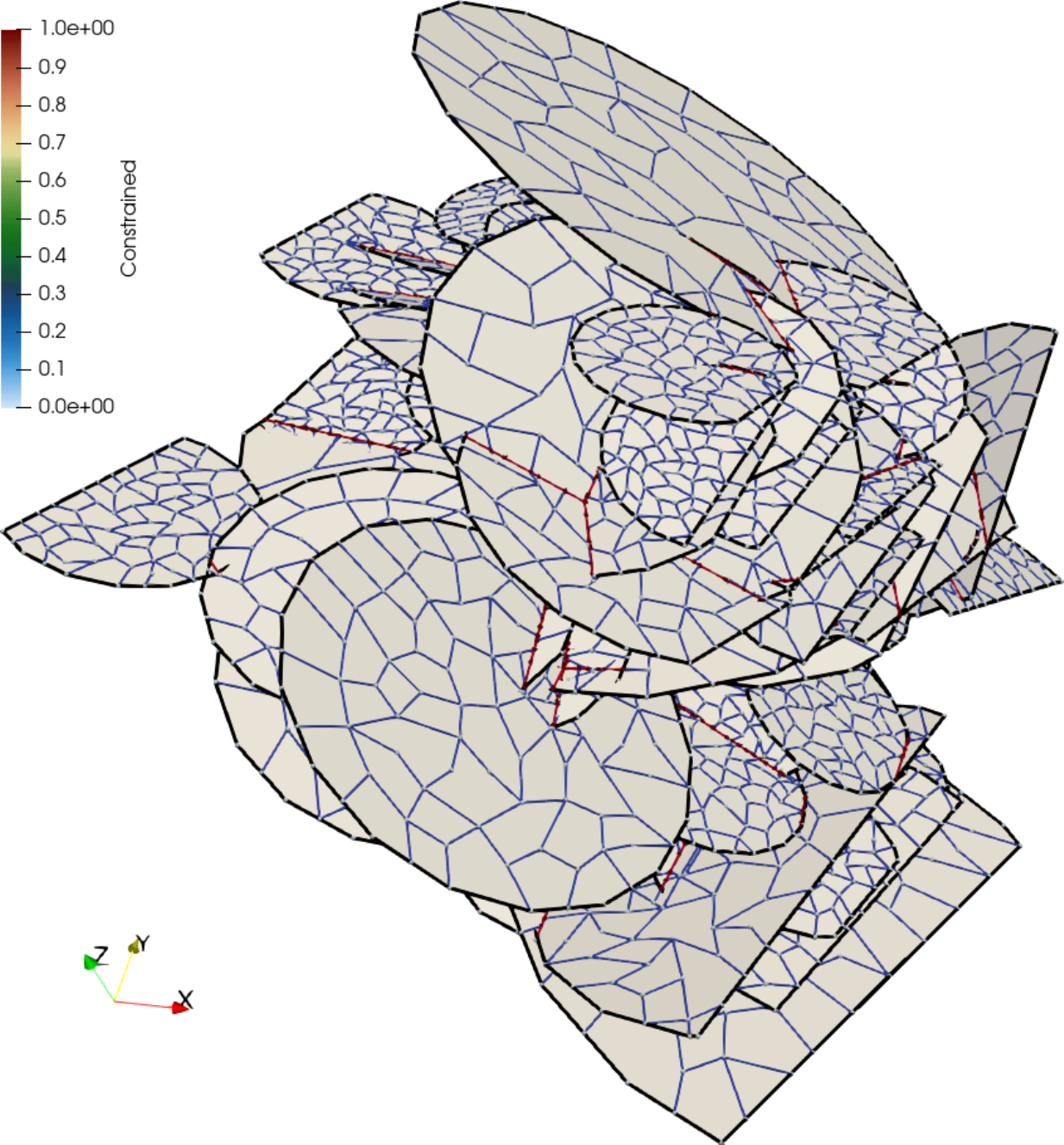}
        \caption{}
        \label{fig:agg:2D:mesh}
    \end{subfigure}
    \begin{subfigure}{0.35\linewidth}
        \includegraphics[width=0.9\linewidth]{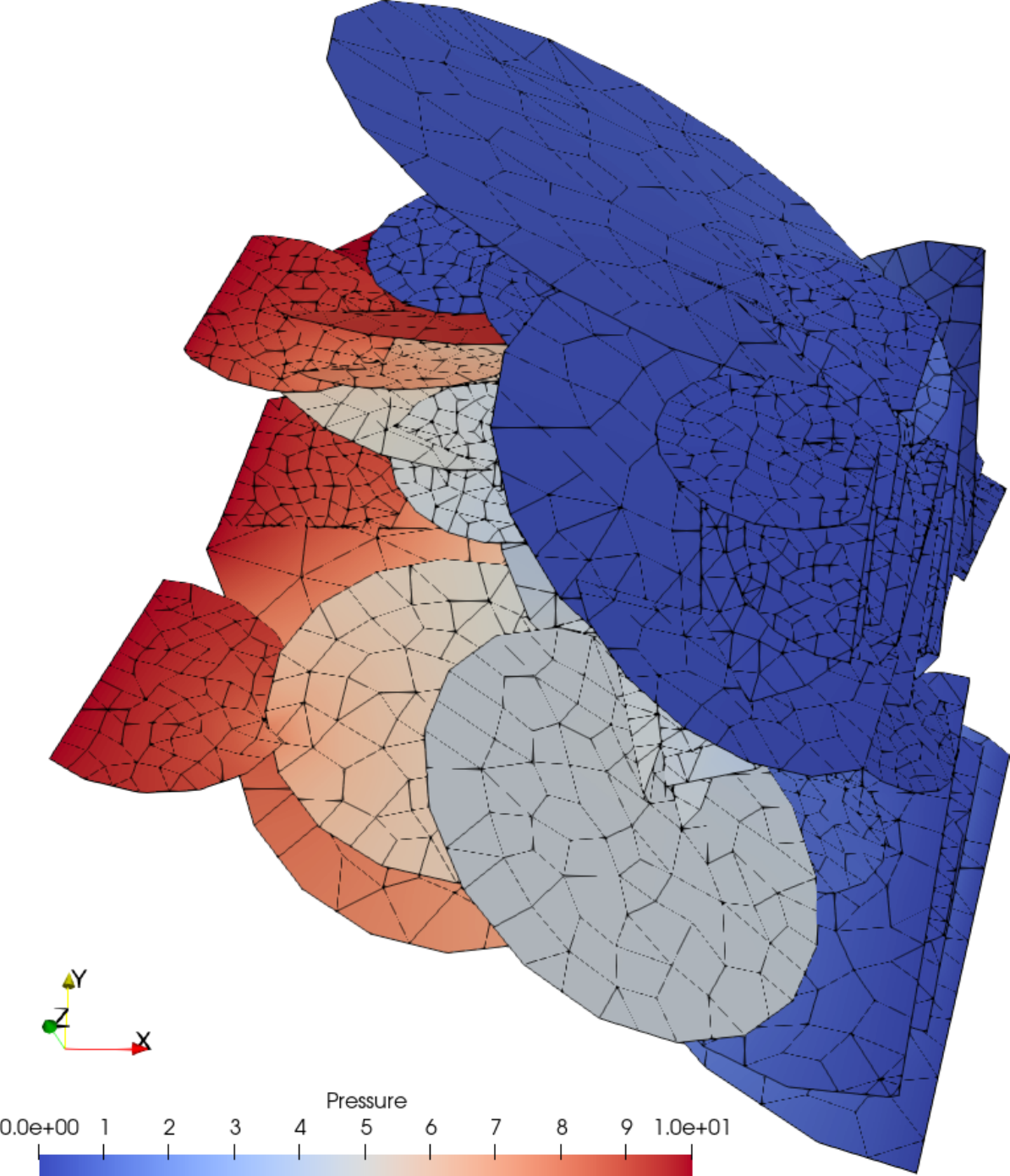}
        \caption{}
        \label{fig:agg:2D:solution}
    \end{subfigure}
    \caption{Test described in \cite{Sorgente2023}. Left: Mesh agglomeration in DFN. Right: DFN coloured by the numerical pressure.}
    \label{fig:agg:2D}
\end{figure}
\begin{figure}[!h]
    \centering
    \begin{subfigure}{0.35\linewidth}
        \includegraphics[width=1.2\linewidth]{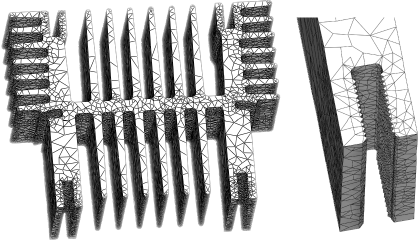}
        \caption{}
    \end{subfigure}\hspace{-27pt}
    \begin{subfigure}{0.35\linewidth}
        \includegraphics[width=0.9\linewidth]{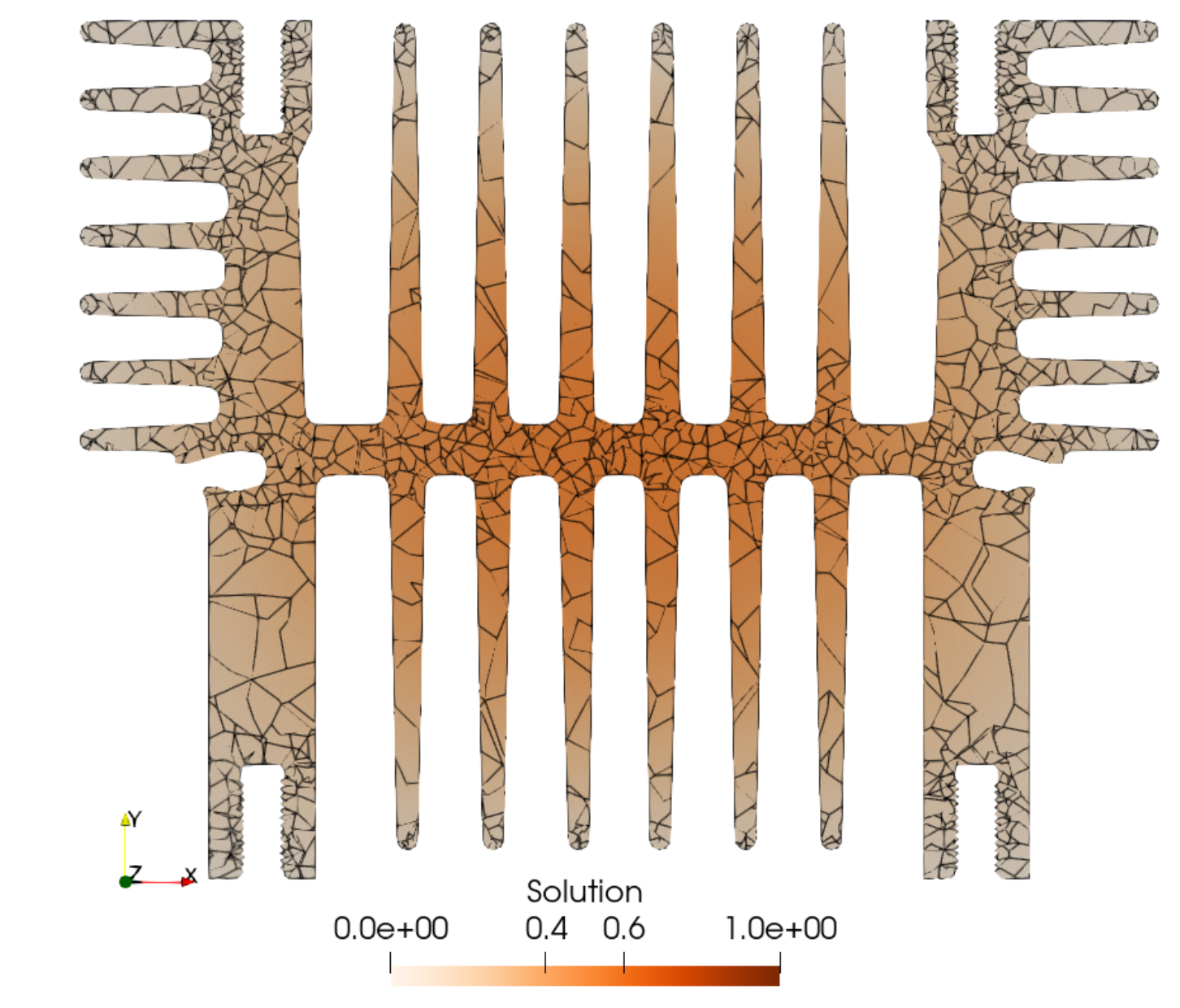}
        \caption{}
    \end{subfigure}
    \caption{Test described in \cite{SORGENTE2025113552}. Left: Mesh agglomeration performed in a mechanical component. Right: Computational domain coloured by the numerical pressure.}
    \label{fig:agg:3D}
\end{figure}
The capability of PolyDiM to generate numerical approximations on generic polytopes has been effectively exploited to enhance tessellation quality through mesh optimization techniques combined with agglomeration strategies \cite{Sorgente2023, SORGENTE2025113552}. In these manuscripts, quality-based optimization strategies to reduce the total number of degrees of freedom associated with a Virtual Element discretization of boundary problems defined over polytopal tessellations are proposed. 

Figures~\ref{fig:agg:2D} and Figures~\ref{fig:agg:3D} illustrate two examples showcasing the application of these techniques to two- and three-dimensional problems, respectively.
In particular, Figure~\ref{fig:agg:2D:mesh} shows the results of this agglomeration-based technique applied to a DFNs, composed by 86 rounded fractures and 159 interfaces, whereas Figure~\ref{fig:agg:2D:solution} illustrates the hydraulic head distribution obtained by solving with Polydim a Darcy problem over this realistic DFN setting.

On the other hand, Figure~\ref{fig:agg:3D} presents the three-dimensional agglomeration-based optimization obtained starting from a tetrahedralization of a complex mechanical component and the corresponding time-dependent solution of a thermal problem defined over this domain obtained with PolyDiM. The reader may refer to \cite{SORGENTE2025113552} for detailed numerical results obtained using this 3D mesh enhancement.

We observe that, in both cases, no relevant differences can be identified between the solution on the non-optimized mesh and the one computed on the agglomerated one, if compared with the huge DOFs reduction obtained.

\subsection{Mesh Refinement}
\begin{figure}[!h]
    \centering
    \begin{subfigure}{0.35\linewidth}
        \includegraphics[width=0.9\linewidth]{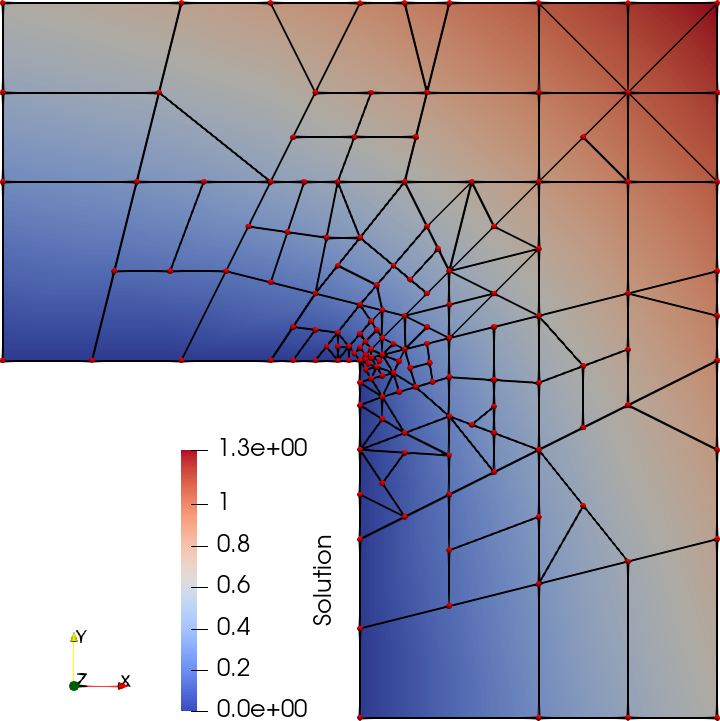}
        \caption{}
    \end{subfigure}
    \begin{subfigure}{0.35\linewidth}
        \includegraphics[width=0.9\linewidth]{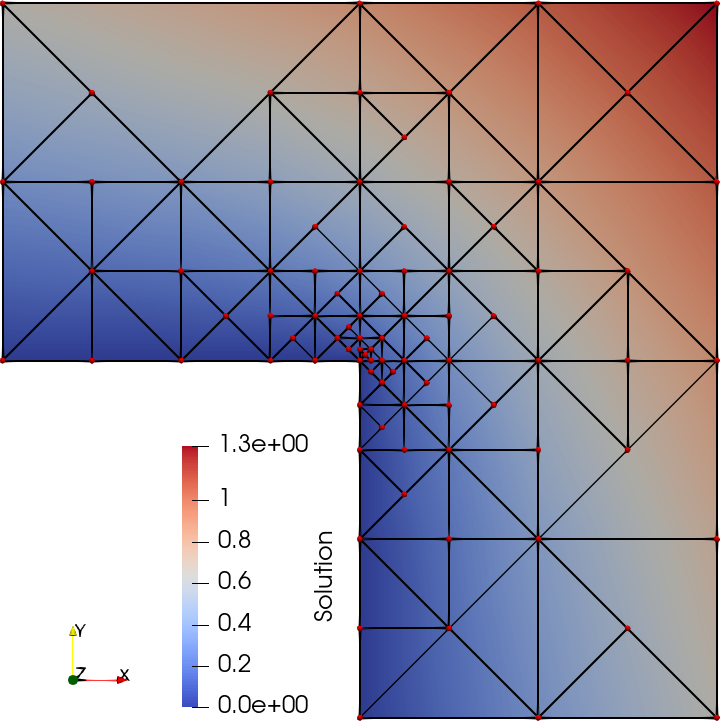}
        \caption{}
    \end{subfigure}
    \caption{Refinement of the L-shaped domain using two different polygonal meshes. The example is detailed in \cite{Vicini2024}.}
    \label{fig:ref:2d}
\end{figure}
\begin{figure}[!h]
    \centering
    \begin{subfigure}{0.35\linewidth}
        \includegraphics[width=1.2\linewidth]{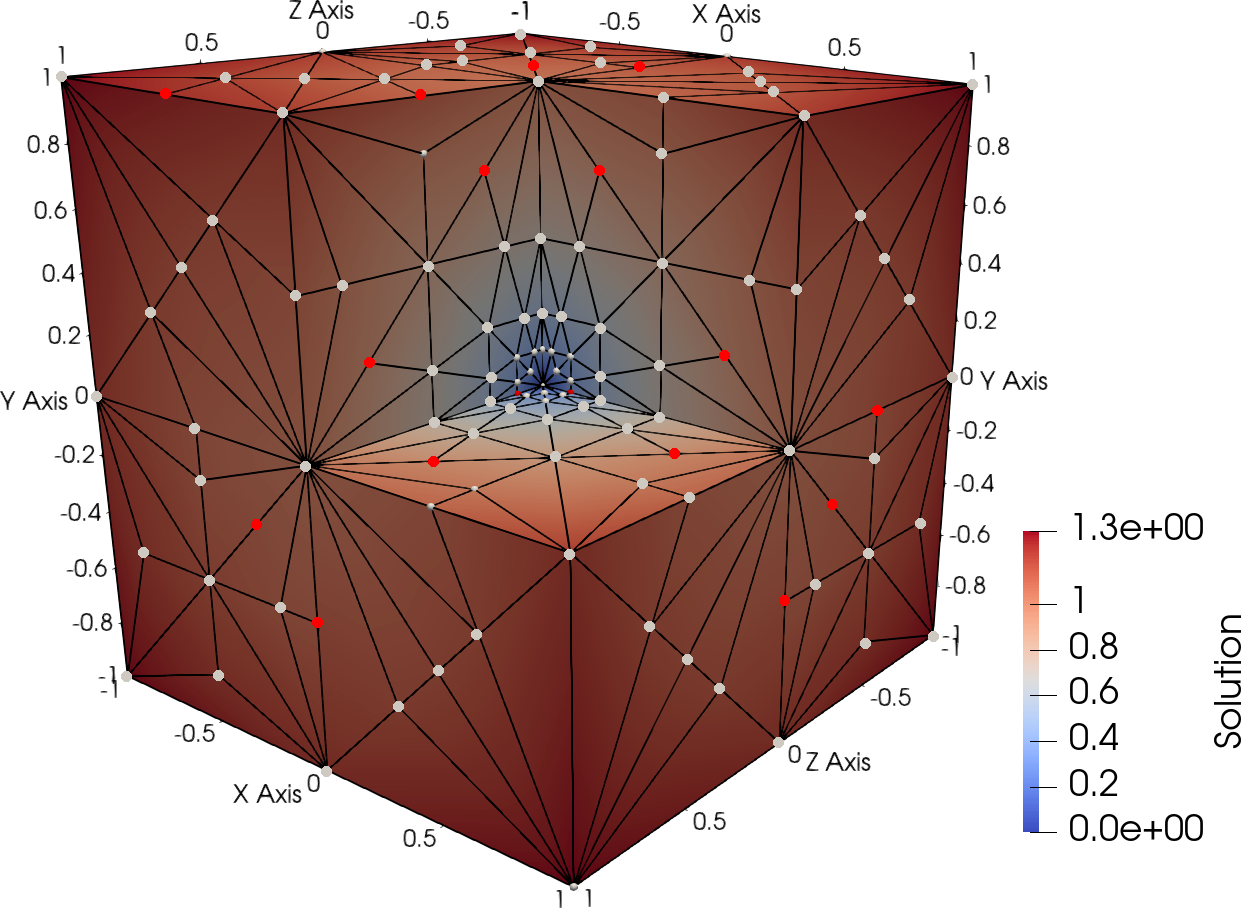}
        \caption{}
    \end{subfigure}\hspace{40pt}
    \begin{subfigure}{0.35\linewidth}
        \includegraphics[width=0.9\linewidth]{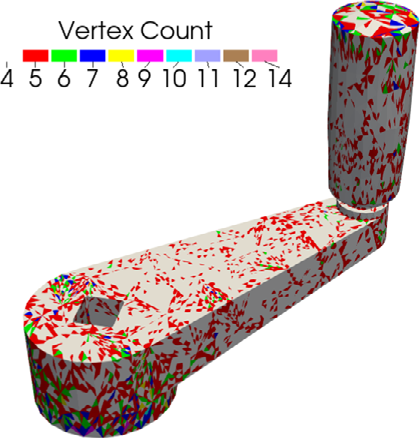}
        \caption{}
    \end{subfigure}
    \caption{These examples are detailed in \cite{Fassino2025}. Left: Refinement of a three-dimensional L-shaped domain. Right: Refinement of a mechanical component. Different colours represent polyhedra with a different number of faces.}
    \label{fig:ref:3d}
\end{figure}

It is well-known that using generic polytopal discrete spaces simplifies the refinement process, as aligned edges or faces can be managed without requiring global mesh updates to maintain conformity. 
The discrete spaces of PolyDiM, combined with the geometric functions provided by GeDiM, support general polygonal and polyhedral mesh refinement on complex geometries. 
Examples of two- and three-dimensional refinements using PolyDiM can be found in \cite{Fassino2023, Vicini2024} and \cite{Fassino2025}, respectively. Specifically, the high order scheme for Adaptive Virtual Element Methods (AVEMs) on triangular meshes has been introduced in \cite{Fassino2023}, whereas work \cite{Vicini2024} details polygonal refinement algorithm to solve elliptic problems with VEM. Finally, the lowest order AVEM scheme on tetrahedral meshes is described by the authors in \cite{Fassino2025}.
Figures~\ref{fig:ref:2d} and \ref{fig:ref:3d} illustrate two cases demonstrating the refinement techniques described in those references. More precisely, Figure~\ref{fig:ref:2d} show the refinement of a two-dimensional L-shaped domain using two different polygonal meshes: the first mesh allows unrestricted aligned edges, while the second constrains the number of aligned edges to a maximum of two.
Figure~\ref{fig:ref:3d}, on the other hand, presents two refined 3D domains: the first is the generalization of the 2D L-shaped problem to the third dimension, while the second is a complex geometry derived from a mechanical component.

\subsection{3D-1D Coupled Problems}
\begin{figure}[!h]
    \centering
    %\hspace{-90pt}
    \begin{subfigure}{0.5\linewidth}
        \includegraphics[width=\linewidth]{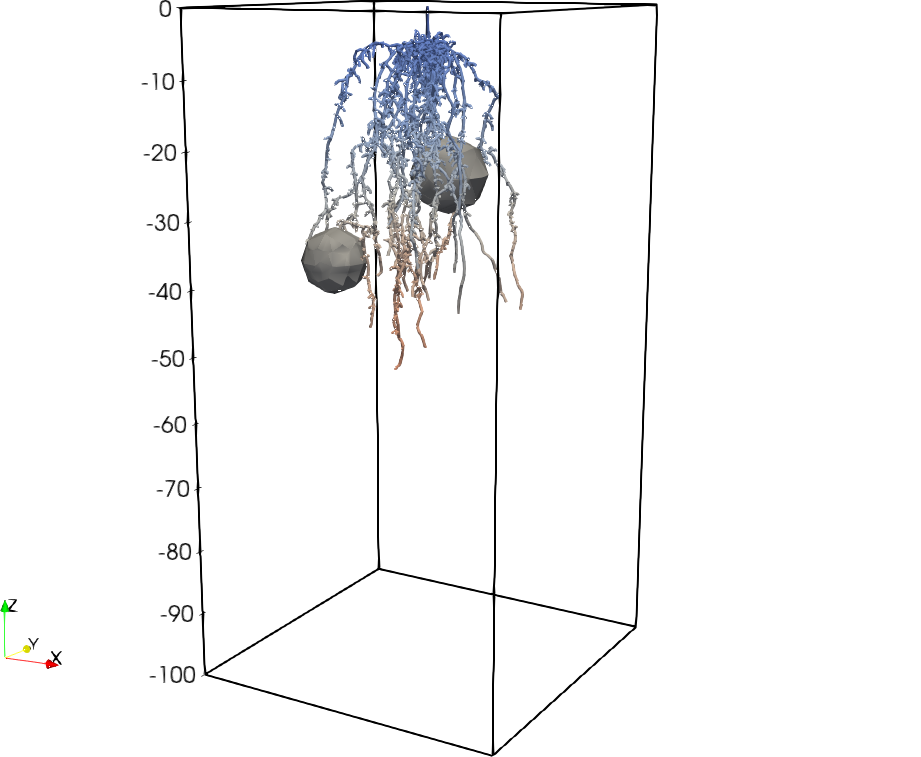}
        \caption{}
    \end{subfigure}\hspace{-20pt}
    \begin{subfigure}{0.5\linewidth}
        \includegraphics[width=\linewidth]{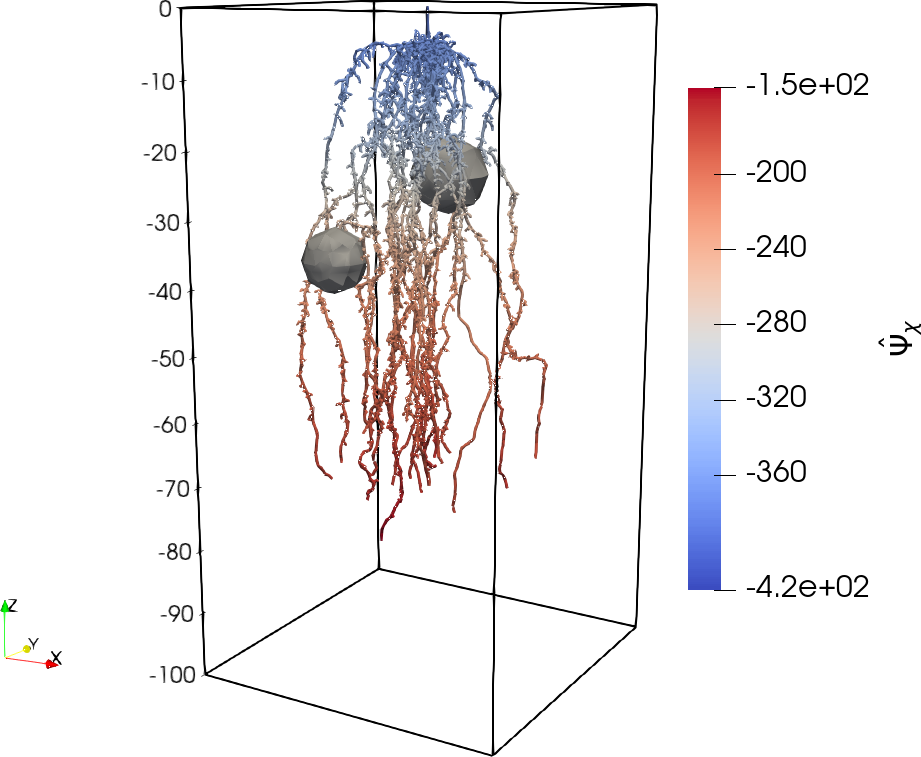}
        \caption{}
    \end{subfigure}
    \caption{The root water pressure across multiple phases of root growth. Left: After 80 days elapsed. Right: 160 days elapsed. This example is detailed in \cite{GrappeinTeora2025}.}
    \label{fig:roots}
\end{figure}
The PolyDiM library has also been successfully applied by the authors in \cite{GrappeinTeora2025} to solve 3D-1D coupled problems. By reducing the 3D-3D coupled problems that models the water uptake by plant roots from soil, described by the Richards equation in the soil domain and by the Stokes equation in the root xylem, the root-soil interactions are simulated on a dynamic root system architecture that develops in a stony soil. Figure \ref{fig:roots} reports a picture illustrating the root water pressure across multiple phases of root growth. The Virtual Element flexibility of handling concave elements and hanging nodes is here exploited by the authors to easily mesh the soil domain near stones and effectively solve the Richards equation in the soil domain.

\section*{Acknowledgements}

The author S.B. and A.B. kindly acknowledges partial financial support provided by PRIN project ``Advanced polyhedral discretisations of heterogeneous PDEs for multiphysics problems'' (No. 20204LN5N5\_003), by PNRR M4C2 project of CN00000013 National Centre for HPC, Big Data and Quantum Computing (HPC) (CUP: E13C22000990001) and the funding by the European Union through project Next Generation EU, M4C2, PRIN 2022 PNRR project P2022BH5CB\_001 ``Polyhedral Galerkin methods for engineering applications to improve disaster risk forecast and management: stabilization-free operator-preserving methods and optimal stabilization methods.''. 
%The author G.T. kindly acknowledges the financial support provided by INdAM-GNCS Project ``Metodi numerici efficienti per problemi accoppiati in sistemi complessi'' (CUP: E53C24001950001). 


\bibliographystyle{IEEEtran}
\bibliography{biblio.bib}

\end{document}